\documentclass[letterpaper, 10 pt, conference]{ieeeconf} 
\IEEEoverridecommandlockouts                             

\usepackage{subfigure}
\usepackage{graphics} 
\usepackage{epsfig} 
\usepackage{scrextend}
\usepackage{times} 
\usepackage{amsmath} 
\usepackage{amssymb} 
\usepackage{color}
\newcommand{\mc}{\mathbb}
\newtheorem{assume}{Assumption}
\newtheorem{lemma}{Lemma}

\newtheorem{remark}{Remark}
\newtheorem{definition}{Definition}

\newtheorem{theorem}{Theorem}
\usepackage{algorithm}
\usepackage{algpseudocode}
\usepackage{Shorthands}
\usepackage{url}
\usepackage{bbold}
\usepackage{mathtools}

 \usepackage{url}
 \usepackage{hyperref}\hypersetup{colorlinks=true, unicode=true, linkcolor=[rgb]{0.10,0.05,0.67}, citecolor=[rgb]{0.10,0.05,0.67}, filecolor=[rgb]{0.10,0.05,0.67}, urlcolor=[rgb]{0.10,0.05,0.67}}
 \usepackage{cite}

\title{\LARGE \bf
Asynchronous Heterogeneous Linear Quadratic Regulator Design
}

\author{Leonardo F. Toso$^\star$, Han Wang$^\star$, and James Anderson% <-this % stops a space
\thanks{$^\star$ Leonardo F. Toso and Han Wang contributed equally to this work. This work is supported in part by  NSF awards 2144634 \& 2231350. The authors are with the Department of Electrical Engineering at Columbia University, New York, NY, 10027, USA. Email: \texttt{\{lt2879, hw2786, james.anderson\}@columbia.edu}.} \thanks{ }% <-this % stops a space
\thanks{
        {\tt\small }}%
\thanks{
        {\tt\small }}%
}

\allowdisplaybreaks

\begin{document}

\maketitle
\thispagestyle{empty}
\pagestyle{empty}

%%%%%%%%%%%%%%%%%%%%%%%%%%%%%%%%%%%%%%%%%%%%%%%%%%%%%%%%%%%%%%%%%%%%%%%%%%%%%%%%
\begin{abstract}
We address the problem of designing an LQR controller in a distributed setting, where $M$ similar but not identical systems share their locally computed policy gradient (PG) estimates with a server that aggregates  the estimates and computes a controller that, on average, performs well on all systems. Learning  in a distributed setting has the potential to offer statistical benefits -- multiple datasets can be leveraged simultaneously to produce more accurate policy gradient estimates. However, the interplay of heterogeneous trajectory data and varying levels of  local computational power introduce bias to the aggregated PG descent direction, and prevents us from fully exploiting the parallelism in the distributed computation. The latter stems from synchronous aggregation, where straggler systems negatively impact the runtime. To address this, we propose an asynchronous policy gradient algorithm for LQR control design. By carefully controlling the ``staleness'' in the asynchronous aggregation, we show that the designed controller converges   
to each system's $\epsilon$-near optimal controller up to a heterogeneity bias. Furthermore, we  prove that our asynchronous approach  obtains exact local convergence at a sub-linear rate. 
\end{abstract}

%%%%%%%%%%%%%%%%%%%%%%%%%%%%%%%%%%%%%%%%%%%%%%%%%%%%%%%%%%%%%%%%%%%%%%%%%%%%%%%%
\vspace{-0.3cm}
\section{Introduction}

In recent years, policy gradient (PG) methods have stood as one of the fundamental pillars underpinning the success of model-free reinforcement learning (RL), offering a versatile framework for learning parameterized policies directly from experience \cite{sutton1999policy, recht2019tour}. Within optimal control, in particular for the linear quadratic regulator (LQR) problem, PG approaches and their \emph{non-asymptotic} performance guarantees have made the task of learning optimal controllers (in a model-free setting), both  possible and systematic. In particular, the authors in \cite{fazel2018global} demonstrated that despite the problem's non-convexity, a derivative-free PG method (i.e., where policy gradients are estimated through simulation data) can converge to the global optimal LQR solution, where the rate of convergence is further proved to be linear in \cite{mohammadi2021convergence}.

Although policy gradient has shown to be effective in learning near-optimal LQR controllers, in a model-free centralized setting \cite{fazel2018global,malik2019derivative,gravell2020learning,mohammadi2019global, toso2023oracle}, it is often assumed that a single agent can access a sufficiently large simulation dataset to produce accurate policy gradient estimates. However,  estimation variance may lead to sub-optimal solutions in a low data regime. To address this, a recent line of work focused on distributed learning for estimation \cite{wang2023fedsysid,toso2023learning,zhang2024sampleefficient,chen2023multi} and control \cite{wang2023model,toso2024meta,zhang2023multi,tang2023zeroth,wang2023fleet}, considers a multi-agent setting, where small contributions (of for example estimated gradients \cite{wang2023model}, system models \cite{wang2023fedsysid}, task weights \cite{zhang2024sampleefficient}, among others) from each agent can be leveraged to achieve statistical benefits. 
In particular, \cite{wang2023model} proposes a federated heterogeneous policy gradient approach to solve the model-free LQR problem. This work shows that by aggregating multiple systems' policy gradient estimates, the sample complexity enjoys a reduction proportional to the number of collaborating systems. However, aggregating heterogeneous PG estimates will inevitably introduce bias to the policy gradient descent direction, as discussed in \cite{wang2023model}.

Critically, the work above implicitly assumes that the policy gradient estimates from all participating agents are promptly available at the beginning of each aggregation step. However, network communication-induced constraints (e.g., rate-limited channels \cite{mitra2024towards}) and agent drop-outs \cite{wang2022fedadmm} may lead to the presence of straggler agents whose late arrival at the server could drastically affect the parallelism of the distributed computation. Namely, an agent with high computational power and/or a fast and reliable communication channel will have its performance throttled as it waits for the idle server (which is waiting for the slower agents) to broadcast the new updated controller. 

To circumvent this limitation, we propose an asynchronous policy gradient approach, where only a batch of the fastest reported estimates at each iteration step are aggregated. This simple modification in the aggregation scheme will mitigate the presence of straggler systems, i.e., by adapting the batch size,  fast agents will not be bothered by long delays incurred  by  waiting for the slow ones to finish their estimates. The trade-off is that the policy gradients produced by the slower agents will be used in the next round of aggregation, they are now out of date, or ``stale''. This staleness may negatively impact the convergence rate of the proposed approach.

Motivated by this, we aim to investigate how the convergence of model-free distributed LQR design is affected by the interplay between stale and heterogeneous PG estimates. In particular, we aim to answer the following questions: Can an asynchronous algorithm produce a controller that is $\epsilon$-near optimal -- even in the presence of staleness and heterogeneous system dynamics? If so, how does the staleness affect the policy gradient convergence? Can we still expect a linear convergence to hold in this setting?

\vspace{-0.2cm}
\subsection{Contributions}

Toward grounding these merits, our main contributions are summarized as follows:

\begin{itemize}
    \item This  is the first work to investigate how aggregating stale PG estimates affects the convergence of an asynchronous model-free distributed LQR design (Algorithm \ref{alg:async_LQR}). We highlight that in our setting, we are dealing with multiple, different systems. We show that the staleness effect can be mitigated by carefully controlling the step size in the PG updates. Moreover, in contrast to a synchronous aggregation scheme, as proposed in \cite{wang2023model}, our approach fully exploits the parallelism in the distributed computation. In particular, it alleviates straggler agents' impact by selecting only a batch of the fastest reported policy gradient estimates at each iteration.   

    \item We first establish the global convergence guarantees of  Algorithm \ref{alg:async_LQR}. To achieve this goal, we derive an upper bound for the staleness term (Lemma \ref{lemma:gradient_bound_interval}), which is approximately of the order of the magnitude of the current iteration step's PG. We prove that our asynchronous algorithm produces a controller that is $\epsilon$-near optimal up to a heterogeneity bias (Corollary \ref{corollary:linear_convergence}). This bias arises due to the heterogeneity among $M$ systems. We demonstrate that staleness impedes global convergence; in particular, the total number of iterations $N$ to achieve each system's $\epsilon$-near optimal controller will be amplified by $\tau^{3/2}_{\max}$ (Corollary ~\ref{corollary:linear_convergence}), where $\tau_{\max}$ is the maximum staleness across systems and PG steps.

    \item We also provide  local convergence guarantees for Algorithm \ref{alg:async_LQR}. Compared to the global convergence analysis, the heterogeneity bias and the staleness effect disappear when converging to a local stationary solution, i.e., our algorithm can exactly locally converge even under the asynchronous aggregation and heterogeneous setting. However, it comes at the cost of a slower convergence, i.e., Algorithm \ref{alg:async_LQR} sub-linearly converges to such fixed point (Corollary \ref{cor:fixed_point_convergence_analysis}). Notably, our tighter local convergence bound demonstrates a linear \emph{speedup} w.r.t. the number of aggregated PG estimates. This improves upon previous work~\cite{nguyen2022federated,toghani2022unbounded}, and thus it is of independent interest to the literature on asynchronous optimization.
\end{itemize}

We validate our theory with numerical experiments, where we demonstrate the benefit of asynchronous aggregation over its synchronous counterpart when dealing with stragglers. Moreover, we show the impact of staleness and heterogeneity on the approach's convergence speed.

\vspace{-0.2cm}
\subsection{Related Work}

\noindent \textbf{Model-free LQR Design:} The setting where a single agent uses its local data to estimate policy gradients and perform controller updates on top of it has been vastly studied to solve the model-free LQR problem \cite{fazel2018global,mohammadi2021convergence,malik2019derivative,gravell2020learning,mohammadi2019global,toso2023oracle}. Although the results on the global and linear convergence are positive and demonstrate the effectiveness of the method, \cite{ziemann2022policy} has shown that PG is very much affected by the limits of control, i.e., poor controlability leads to arbitrarily noisy gradient estimates. On the other hand, \cite{wang2023model,toso2024meta,tang2023zeroth} have demonstrated the value of collaboration, in a multi-agent setting, to reduce the variance in the estimated gradient and achieve sample efficiency when learning LQR controllers. Most relevant to our work is \cite{wang2023model}, where the authors consider a synchronous policy gradient approach to tackle the model-free LQR problem. In contrast to \cite{wang2023model}, our work considers an asynchronous aggregation scheme where the impact of straggler agents is mitigated and the effect of aggregating stale PG estimates is thoroughly characterized in our local and global convergence analysis.

\vspace{0.2cm}
\noindent \textbf{Asynchronous Optimization:} Asynchronous stochastic gradient descent (SGD) has been a topic of study in stochastic optimization over the past decade, where many papers \cite{agarwal2011distributed,chaturapruek2015asynchronous,stich2021critical} investigate the connection of large batches and staleness in the ergodic convergence rate of such approach. Most relevant to our work are the recent papers on asynchronous distributed learning \cite{nguyen2022federated,toghani2022unbounded,karimi2016linear,ma2023csmaafl,fabbro2024dasa,adibi2024stochastic}, where in contrast to asynchronous SGD, the heterogeneity in the local data distribution and local data privacy concerns are taking into account. In contrast to these previous work, our paper consider an asynchronous policy gradient approach to solve the LQR optimal control problem. In such setting, not only the convergence guarantees are importantly characterized, but the per-iteration closed-loop stability of the collaborating agents under the designed controller is also required. Moreover, our theoretical guarantees reveals a linear speedup with respect to the number of aggregated policy gradient estimates in the local convergence rate, which is of independent interest to provide tighter bounds in \cite{nguyen2022federated,toghani2022unbounded}.

\vspace{-0.1cm}
\subsection{Notation}

Let $[M]$ denote the set of integers $\{1,\hdots, M\}$. We use $\mathcal{J}(K)$ to denote the LQR cost for an arbitrary system with problem parameters $(A,B,Q,R)$. When required,  $\mathcal{J}^{(i)}(K)$ denotes the LQR cost for a specific tuple $(A^{(i)}, B^{(i)}, Q^{(i)}, R^{(i)})$, where $i \in [M]$. The spectral radius of a square matrix is $\rho(\cdot)$, and $\sigma_{\min}(\cdot)$ denotes the minimum singular value. Unless otherwise stated, $\|\cdot\|$ is the spectral norm. We use  $\mathcal{O}(\cdot)$ to omit constant factors in the argument.

\section{Problem Formulation}\label{sec:problem_formulation}

Let us begin with the standard setup of multi-agent LQR design for heterogeneous systems \cite{wang2023model,toso2024meta}. Consider $M$ discrete-time and linear time-invariant (LTI) dynamical systems over an infinite time-horizon, described by 
\begin{align}\label{eq:LTI_systems}
    x^{(i)}_{t+1} = A^{(i)}x^{(i)}_t + B^{(i)}u^{(i)}_t, \quad  \forall i \in [M], 
\end{align}
where $A^{(i)} \in \mathbb{R}^{n_x\times n_x}$, $B^{(i)} \in \mathbb{R}^{n_x\times n_u}$ are the system matrices, and $x^{(i)}_t \in \mathbb{R}^{n_x}$, $u^{(i)}_t \in \mathbb{R}^{n_u}$ denote the state and control input of system $i$ at time instant $t$, respectively. The initial state $x^{(i)}_0$ of \eqref{eq:LTI_systems} is drawn from an arbitrary distribution $\mathcal{X}_0$ that satisfies Assumption \ref{assumption:initial_state_distribution} below.  In this work we specifically account for the fact that the $M$ systems are \emph{not identical}, i.e., in general $A^{(i)}\neq A^{(j)}$ and $B^{(i)}\neq B^{(j)}$. We quantify the level of heterogeneity at the end of this section.

From the perspective of the $i^{\mathrm{th}}$ system, the goal is to design an optimal static state feedback controller $K^\star_i \in \mathcal{K}^{(i)} := \{K \in \mathbb{R}^{n_u\times n_x} \mid \rho(A^{(i)} - B^{(i)}K) < 1\}$, that provides a control policy $u^{(i)}_t = -K^\star_ix^{(i)}_t$ that, subject to~\eqref{eq:LTI_systems}, minimizes the quadratic cost function
\begin{align} \label{eq:LQR_cost}
 \mathcal{J}^{(i)}(K) := {\mc{E}}\left[\sum_{t=0}^{\infty} x^{(i)\top}_t \left(Q^{(i)} + K^\top R^{(i)} K\right) x^{(i)}_t\right],
\end{align}
where $Q^{(i)} \in \mathbb{R}^{n_x \times n_x}$ and $R^{(i)} \in \mathbb{R}^{n_u\times n_u}$ denote the cost matrices, and the expectation is with respect to $x^{(i)}_0 \sim \mathcal{X}_0$. 

\begin{assume}[Initial state distribution]\label{assumption:initial_state_distribution}
The initial state distribution $\mathcal{X}_0$ satisfies $\mc{E}[x_0^{(i)}]=0$ (i.e., zero mean) with covariance $\Sigma_0 = \mc{E}[x_0^{(i)}x_0^{(i)\top}] \succ \mu I_{n_x}$  for some $\mu > 0$.\footnote{Assumption \ref{assumption:initial_state_distribution} implies that the covariance matrix has full rank. It is a standard assumption in PG-LQR literature \cite{fazel2018global, malik2019derivative,gravell2020learning,mohammadi2019global} and guarantees that all stationary solutions are global optima.}
\end{assume}

This work considers the model-free setting where the tuple $(A^{(i)}, B^{(i)}, Q^{(i)}, R^{(i)})$ is \emph{unknown} and each system only has limited access to simulation data. Thus, designing $K^\star_i$ through the well-established Riccati equation \cite{hewer1971iterative} is not possible. As in \cite{fazel2018global, malik2019derivative}, we must resort to derivative-free policy gradient approaches to optimize \eqref{eq:LQR_cost}.

However, due to limited local trajectory data, accurate policy gradient estimations in the single-agent setting may require more data than is available. As proposed in \cite{wang2023model}, instead of designing $K^\star_i$ for each system $i \in [M]$, we consider the problem of leveraging  simulation data from multiple (differing but ``similar'') systems in order to compute a controller $\Kbar^\star \in \mathcal{K} :=  \cap \mathcal{K}^{(i)}$ that, i) stabilizes each system, ii) on average, performs well for all of them. Moreover, such an approach should be sample efficient, in the sense that a small data contribution from multiple systems can be leveraged into a larger dataset to perform more accurate PG estimates. With regards to ii), $\Kbar^\star$ minimizes
\begin{align} \label{eq:avg_cost}
\Kbar^\star :=  & \argmin_{K \in {\mathcal{K}}} \left\{ \bar{\mathcal{J}}(K) := \frac{1}{M}\sum_{i=1}^M \mathcal{J}^{(i)}\left(K \right) \right\},
\end{align}
subject to \eqref{eq:LTI_systems}. To solve \eqref{eq:avg_cost}, we first consider an arbitrary initial stabilizing controller $\Kbar \in \mathcal{K}$ and step-size $\eta \in \mathbb{R}_{>0}$, and iteratively perform policy gradient updates of the form: 
\[
\Kbar \leftarrow \Kbar - \eta\widehat{\nabla}\bar{\mathcal{J}}(\Kbar),
\] 
where $\widehat{\nabla}\bar{\mathcal{J}}(\cdot)$ is an estimate of the true gradient ${\nabla}\bar{\mathcal{J}}(\cdot)$. 

In \cite{wang2023model}, at each iteration $n\in \{0,1,2,\ldots\}$, all $M$ policy gradient estimates generated by the systems, $\widehat{\nabla}\mathcal{J}^{(i)}(\Kbar_n)$, are aggregated to produce $\widehat{\nabla}\bar{\mathcal{J}}(\Kbar_n)$. This is then used to compute the new policy $\Kbar_{n+1}$. It is implicitly assumed that the policy gradients are all available to produce $\widehat{\nabla}\bar{\mathcal{J}}(\Kbar_n)$. Such an assumption (i.e., synchronous aggregation) does not take into account network communication effects. While synchronously aggregating such estimates may offer sample efficiency \cite{wang2023model}; however, the presence of straggler systems will prevent  full exploitation of the parallelism in the  distributed computation. To circumvent this limitation, we consider an asynchronous distributed LQR design, where, at each iteration step $n$, only a subset $[b_s] \subseteq [M]$ of the first reported PG estimates contribute to the update of $\Kbar_{n+1}$, i.e., 
\vspace{-0.2cm}
\begin{align}\label{eq:async_update}
    \Kbar_{n+1} = \Kbar_{n} - \frac{\eta}{b_s}\sum_{s = 1}^{b_s}\widehat{\nabla}\mathcal{J}^{(s)}(\Kbar_{n - \tau_s(n)}),
\end{align}
where $\tau_s(n)$ denotes the staleness in the controller that system $s \in [b_s]$ possesses when locally estimating its policy gradient at  step $n$. Due to heterogeneity, in the synchronous setting \cite{wang2023model}, $\underset{n \to \infty}{\limsup} \Kbar_{n+1} \neq \Kbar^\star_i$, and $\mathcal{J}^{(i)}(\Kbar_n) - \mathcal{J}^{(i)}(K^\star_i)$ is characterized by the heterogeneity level across systems.

\begin{figure}
    \centering
    \includegraphics[width=0.45\textwidth]{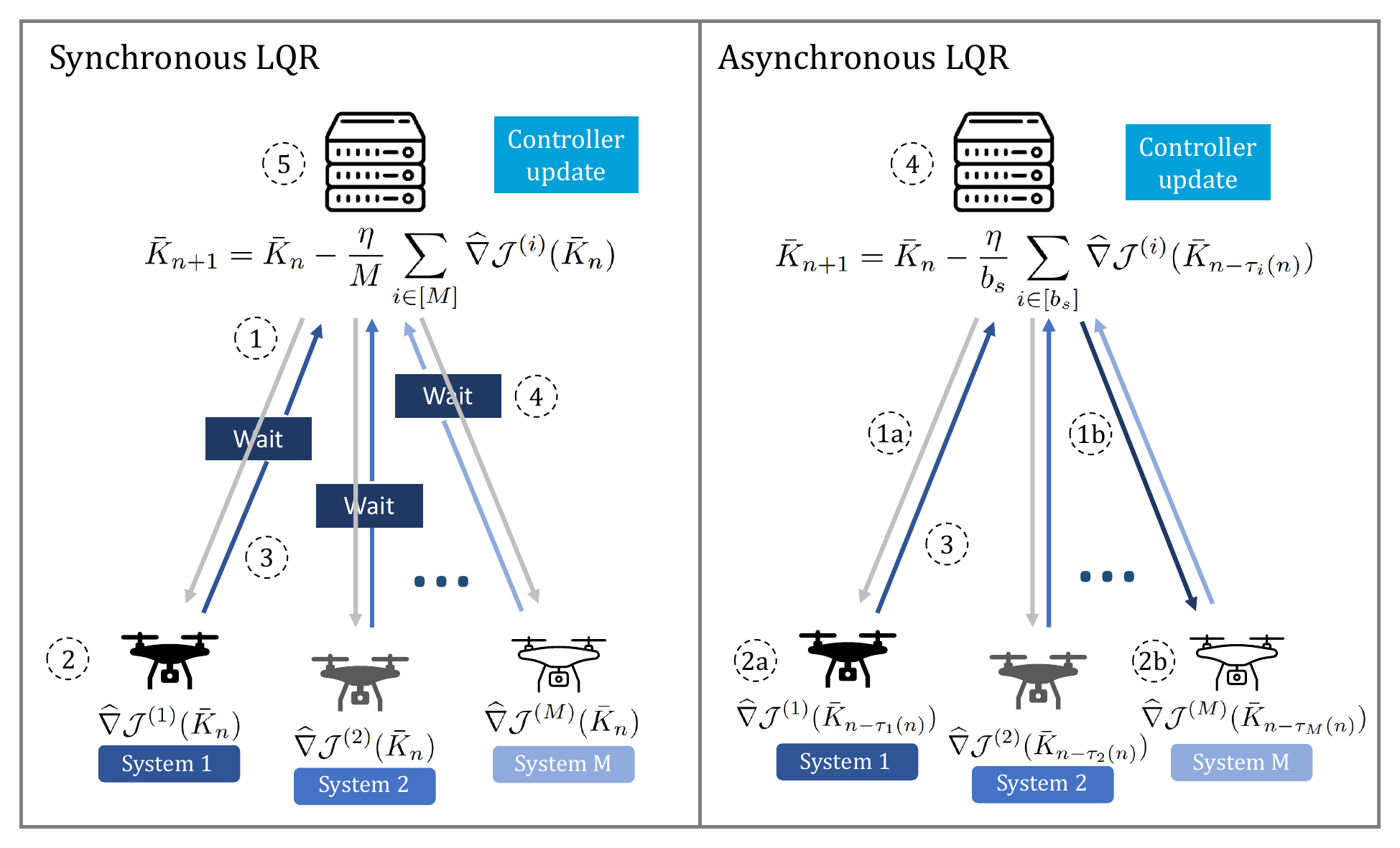}
    \caption{Illustrative schematic to compare the synchronous and asynchronous PG approaches for the model-free LQR problem. Figure inspired by \cite{ma2023csmaafl}.} 
    \label{fig:sync_async_diagram}
\end{figure}

Intuitively, besides heterogeneity, aggregating over stale estimates may also produce sub-optimal solutions to the LQR design. Therefore, in this work, we analyze the interplay between heterogeneity and the staleness in the convergence of the policy update \eqref{eq:async_update}. Figure \ref{fig:sync_async_diagram} illustrates the comparison between the synchronous and asynchronous PG for the model-free LQR problem. As we can see on the left-hand side, after the server communicates the updated controller (step 1), it needs to wait (step 4), for all systems to complete their estimates (step 2) and communicate back to the server (step 3), before performing the next controller update (step 5). On the other hand, in the asynchronous PG, despite of the staleness in the updated controller (step 1a or 1b) that the system has to estimate its gradient (step 2a or 2b), only the first $b_s$ PG estimates reported back to the server (step 3) are aggregated in the controller update (step 4). The impact of straggler systems is then alleviated by adapting $b_s$. 

Before presenting the algorithm that implements such asynchronous policy gradient update, we first need to define the sub-level set of stabilizing controllers $\mathcal{S} \subseteq \mathcal{K}$, as well as formally state the aforementioned requirement of initially iterating from a stabilizing controller.

\begin{definition}\label{def:stabilizing_set} Let $\Kbar_0$ and $K^\star_i$ be the initial and optimal stabilizing controllers of system $i \in [M]$, respectively. The stabilizing sub-level set of $\mathcal{K}$ is $\mathcal{S} \triangleq \cap \mathcal{S}^{(i)}$, with
    \begin{align*}
        \mathcal{S}^{(i)}:= \left\{K\; | \; \mathcal{J}^{(i)}(K) - \mathcal{J}^{(i)}(K^\star_i) \leq \gamma \Delta^{(i)}_0\right\},  
    \end{align*}
where $\Delta^{(i)}_0 \triangleq \mathcal{J}^{(i)}(\Kbar_0) - \mathcal{J}^{(i)}(K^\star_i)$ denotes the initial distance to optimality, and $\gamma \geq 1$ is an arbitrary scalar.
    
\end{definition}

\begin{assume}[Initial stabilizing controller $\Kbar_0$] \label{assumption:initial_controller}The initial controller $\bar{K}_0$ stabilizes all systems, i.e. $\bar{K}_0 \in \mathcal{S}$.
\end{assume}

The assumption on the initial stabilizing controller is standard in policy gradient methods for LQR design \cite{fazel2018global, gravell2020learning,mohammadi2019global}. If $\Kbar_0$ is not stabilizing for all the systems \eqref{eq:LTI_systems}, then \eqref{eq:async_update} will not produce a stabilizing controller, since $\mathcal{J}^{(i)}(\Kbar_0)$ is undefined for the corresponding unstabilized systems. In addition, we emphasize that although $\Kbar_0$ stabilizes \eqref{eq:LTI_systems} $\forall i$, it may provide a sub-optimal cost $\bar{\mathcal{J}}(\Kbar_0) \geq \bar{\mathcal{J}}(\Kbar^\star)$.

\section{Asynchronous Policy Gradient Algorithm}

The proposed asynchronous policy gradient algorithm, as described in Algorithm \ref{alg:async_LQR}, implements the controller update in \eqref{eq:async_update}. The subset $[b_s] \subseteq [M]$ of policy gradient estimates are computed through a two-point zeroth-order gradient estimation approach, as detailed in Algorithm \ref{alg:ZO2P}. Upon initializing all systems $i \in [M]$ with an initial stabilizing controller $\Kbar_0$, i.e., step 2 of Algorithm \ref{alg:async_LQR}, in parallel, each system estimate its policy gradient using its local simulation data in step 3. The zeroth-order estimation in Algorithm \ref{alg:ZO2P} performs an empirical estimation of the first-order Gaussian Stein's identity \cite{stein1972bound}. That is, given a smoothing radius $r$, the current controller $K$ is perturbed by a random matrix $U$, such that $\|U\|_F = r$, to produce $K^1 = K + U$ and $K^2 = K - U$. Such smoothing controllers are then played by the $i^{\mathrm{th}}$ system to collect simulation data and compute the costs $\mathcal{J}^{(i)}(K^1)$, and $\mathcal{J}^{(i)}(K^2)$. Then, by averaging over $m$ samples, Algorithm \ref{alg:ZO2P} returns a biased empirical estimation $\widehat{\nabla}\mathcal{J}(K)$ of $\mc{E}\left[\nabla \mathcal{J}^{(i)}(K)\right] = \mc{E}\left[\frac{n_xn_u}{2r^2}(\mathcal{J}^{(i)}(K^1) - \mathcal{J}^{(i)}(K^2))U\right]$.

\begin{algorithm}
\caption{Asynchronous Distributed LQR} 
\label{alg:async_LQR}
\begin{algorithmic}[1]
\State \textbf{Input:} Stabilizing controller $\Kbar_0$, step-size $\eta$, batch size $b_s$, iterations $N$, smoothing radius $r$, and number of samples $m$.\
\State \textbf{Initialize} the local controllers $K_i = \bar{K}_0$ $\forall i \in [M]$, batch and iteration counters $s = n = 0$, and $\overline{\nabla} \leftarrow 0$.
\State \textbf{In parallel} compute and send $\widehat{\nabla}_i = \texttt{ZO}(K_i,r,m)$ to the server $\forall i \in [M]$.\ 
\State \textbf{While} $n < N$\
\State \quad \textbf{If} the server receives an estimate \textbf{then}\
\State \quad \quad Accumulate $\overline{\nabla} = \overline{\nabla} + \widehat{\nabla}_i$,  $s \leftarrow s +1$,
\State \quad \textbf{If} $s = b_s$ \textbf{then}
\State \quad \quad  $\bar{K}_{n+1} = \bar{K}_n - \frac{\eta}{b_s}{\overline{\nabla}}$,\;\  $n \leftarrow n + 1$, $\overline{\nabla} \leftarrow 0$, $s \leftarrow 0$,\
\State \quad \quad \textbf{If} system $i \in [M]$ is \textbf{done}, \textbf{then} $K_i \leftarrow \bar{K}_{n+1}$ \textbf{and}
\State \quad \quad   \quad Compute $\widehat{\nabla}_i = \texttt{ZO}(K_i,r,m)$,\
\State \quad \quad \quad Send $\widehat{\nabla}_i$ to the server,\
\State \textbf{Output:} $\bar{K}_N$.
\end{algorithmic}
\end{algorithm}

Once a system $i \in [M]$ is done estimating its PG, it sends the estimate to a server that accumulates them (step 6). A policy gradient update \eqref{eq:async_update} is only performed when the number of accumulated gradient estimates is equivalent to the batch size $b_s$ (step 9), indicating that the $b_s$ fastest reported PG estimates at iteration $n$ are ready to be aggregated. After $N$ iterations of steps 4-11, Algorithm \ref{alg:async_LQR} returns $\Kbar_N$. In Section \ref{sec:theoretical_guarantees}, we characterize the properties of $\Kbar_N$ based on a local and global convergence rates. We re-emphasize that Algorithm \ref{alg:async_LQR} aggregates stale policy gradient estimates in step 8. In addition, due to the heterogeneity, such staleness in the controller that each system access to perform its gradient estimate is later assumed to be bounded\footnote{This can be relaxed to unbounded staleness in the homogeneous asynchronous distributed learning setting, as discussed in \cite{mishchenko2022asynchronous}.} in Assumption \ref{assumption:bounded_staleness}.

\begin{remark} We exploit a two-point gradient estimation approach since it offers a lower estimation variance compared to the one-point counterpart \cite{malik2019derivative}. Moreover, for simplicity and by following \cite{malik2019derivative,toso2023oracle},  in Algorithm \ref{alg:ZO2P}, it is implicitly assumed to have access to the true infinite horizon costs $\mathcal{J}^{(i)}(K^1)$ and $\mathcal{J}^{(i)}(K^2)$. However, since such quantities are lower bounded by any finite-horizon approximation, our results can be readily extended to that setting as well \cite{gravell2020learning}.  \end{remark}

\begin{algorithm}
\caption{\texttt{ZO}: Two-point Zeroth-order Estimation} 
\label{alg:ZO2P}
\begin{algorithmic}[1]
\State \textbf{Input:} Stabilizing controller $K$, number of samples $m$ and smoothing radius $r$.
\State \textbf{for} all samples $l \in [m]$ \textbf{ do } \
\State \quad \textbf{Draw} $U_l \in \mathbb{R}^{n_u\times n_x}$, such that $\|U_l\|_F = r$,
\State \quad \textbf{Smooth} controllers: $K^1_l=K+U_l$ and $K^2_l=K-U_l$,
\State \quad \textbf{Compute} and \textbf{store} costs $\mathcal{J}(K^1_l)$, and $\mathcal{J}(K^2_l)$,
\State \textbf{end for}\
\State \textbf{Return} $\widehat{\nabla}\mathcal{J}(K)= \frac{n_x n_u}{2 r^2 m} \sum_{l=1}^{m}  (\mathcal{J}(K^1_l) - \mathcal{J}(K^2_l)) U_l$
\end{algorithmic}
\end{algorithm}

\begin{assume}[Bounded staleness]\label{assumption:bounded_staleness}  For any system $i \in [M]$ and iteration $n$, $\tau_i(n) \leq \tau_{\max}$, for some $\tau_{\max} \in [1,\infty)$.
\end{assume}

The above assumption is common in the convergence analysis of asynchronous stochastic gradient descent algorithms \cite{ nguyen2022federated, toghani2022unbounded, koloskova2022sharper}. It guarantees that the asynchronous aggregation in Algorithm \ref{alg:async_LQR} is performed within a finite time. 

Before jumping to the convergence and stability analysis of Algorithm \ref{alg:async_LQR}, we first revisit some properties of the policy gradient LQR \cite{fazel2018global} and heterogeneity bound \cite{wang2023model} that are instrumental in deriving the main results of this work.

\begin{lemma}[Local smoothness]\label{lemma:smoothness} Given a pair of stabilizing controllers $K, K^{\prime} \in \mathcal{S}$ such that $\|K^{\prime} -K\|_F <\infty$,  the LQR gradient is $\barhgrad$-Lipschitz, i.e., 
\begin{align*}
&\left\|\nabla \mathcal{\mathcal{J}}^{(i)}\left(K^{\prime}\right)-\nabla \mathcal{J}^{(i)}(K)\right\|_F \leq h_{\text {grad}}\|K^{\prime} -K\|_F,
\end{align*}
where $h_{\text {grad}}$ depends on the LQR problem parameters.
\end{lemma}

The proof of the LQR gradient's local smoothness was first introduced in \cite{fazel2018global}, and the explicit expression of $h_{\text {grad}}$ was further provided in \cite[Appendix D.1]{wang2023model}.

\begin{lemma}[Gradient dominance]\label{lemma:gradient_dominance}  Let $K^\star_i$ be the LQR optimal controller associated with system $i \in [M]$. Given a stabilizing controller $K \in \mathcal{S}$ the squared norm of the LQR gradient is lower bounded as follows:
\begin{align*}
      \|\nabla \mathcal{J}^{(i)}(K)\|_F^2 \geq \lambda\left(\mathcal{J}^{(i)}(K)-\mathcal{J}^{(i)}\left(K^\star_i\right)\right), \;\ \forall i \in [M],
\end{align*}
where $\lambda = 4\mu^2 \underset{i \in [M]}{\max} \sigma_{\min }(R^{(i)})/\|\Sigma_{K^\star_i}\|$ denotes the gradient dominance constant and $\Sigma_{K^\star_i} = \mc{E}[x_t^{(i)}x_t^{(i)\top}]$ corresponds to the state covariance matrix incurred by playing  \eqref{eq:LTI_systems} with its corresponding optimal controller.
\end{lemma}

Lemmas \ref{lemma:smoothness}-\ref{lemma:gradient_dominance} are paramount to prove the global convergence of the policy gradient LQR in the single-agent setting \cite{fazel2018global}. Even though, the gradient dominance property of each ${\mathcal{J}}^{(i)}(K)$ does not imply the same property to the the average cost $\bar{\mathcal{J}}(K)$, i.e., due system and cost heterogeneity. We can still leverage such result when characterizing the distance to optimality, i.e., $\Delta^{(i)}_N := \mathcal{J}^{(i)}(\Kbar_N) - \mathcal{J}^{(i)}(K^\star_i)$, for all systems $i \in [M]$, in Section \ref{sec:theoretical_guarantees}.

We now quantify the level of heterogeneity between the $M$ systems. Note that we could include a common bound on the spectral norm difference among system and cost matrices at the expense of less precise downstream results.
\begin{assume}\label{assumption:bounded_het}  There exist positive scalars $\epsilon_A, \epsilon_B, \epsilon_Q, \epsilon_R$, such that system and cost heterogeneity is bounded. That is,
\begin{align*}
&\underset{i\neq j}{\max} \lVert  A^{(i)} -A^{(j)}\rVert \leq \epsilon_A,\;\ 
\underset{i\neq j}{\max}\lVert B^{(i)} -B^{(j)} \rVert\leq \epsilon_B, \\
&\underset{i\neq j}{\max} \lVert  Q^{(i)} - Q^{(j)}\rVert \leq \epsilon_Q, \;\
\underset{i\neq j}{\max}\lVert R^{(i)} -R^{(j)} \rVert\leq \epsilon_R.
\end{align*}

\end{assume}

\begin{lemma}[Gradient heterogeneity \cite{toso2024meta}]\label{lemma:gradient_heterogeneity}
Given a pair of distinct systems $i\neq j \in [M]$, following the dynamics in \eqref{eq:LTI_systems}, and a stabilizing controller $K \in \mathcal{S}$. The following holds, 
\begin{align*}
    \|\nabla \mathcal{J}^{(i)}(K) - \nabla \mathcal{J}^{(j)}(K)\|^2_F \leq {\epsilon}_{\text{het}}.
\end{align*}
where  ${\epsilon}_{\text{het}} $ scales quadratically with the system and cost heterogeneity levels $ (\epsilon_A, \epsilon_B, \epsilon_Q, \epsilon_R)$ of Assumption \ref{assumption:bounded_het}. 
\end{lemma}

The proof of the above lemma as well as the explicit expression of ${\epsilon}_{\text{het}}$ is presented in \cite[Appendix 6.2]{toso2024meta}. Note that, for a small gradient heterogeneity level ${\epsilon}_{\text{het}}$, this lemma conveys that the PG descent directions of systems $i \neq j \in [M]$ are close, which is then also close to the descent direction of the average cost. This lemma is crucial to quantify the effect of heterogeneity in the convergence and stability analysis of our asynchronous policy gradient aggregation.

\section{Convergence and Stability Analysis}\label{sec:theoretical_guarantees}
We now present the main theoretical results of this work. First, we show that Algorithm \ref{alg:async_LQR} can exactly converge to the local optimum at a sub-linear convergence rate. Second, we provide global convergence guarantees for our proposed approach. In the global convergence analysis, due to heterogeneity, our algorithm will converge to a ball that contains each system's optimal controller. The size of the convergence ball depends on the heterogeneity level among systems. We demonstrate that, even in the presence of a staleness, this convergence is achieved at a linear rate with respect to the tolerance level. Moreover, we establish a linear convergence rate that has a dependence on the maximum staleness $\tau_{\max}$.

\subsection{Local convergence guarantees}

First, the local convergence of Algorithm \ref{alg:async_LQR} is characterized through the ergotic convergence rate, i.e., how $\frac{1}{N}\sum_{n=0}^{N-1}\mc{E}\|\nabla \bar{\mathcal{J}}(\Kbar_n)\|^2_F$ scales with the number of iterations $N$, batch size $b_s$, heterogeneity $\het$ and staleness $\tau_{\max}$.

\begin{theorem}\label{thm:fixed_point_convergence_analysis} 
Let Assumptions \ref{assumption:initial_state_distribution}-\ref{assumption:bounded_het} hold. Suppose the step-size satisfies $\eta \leq \barhgrad\eta_{\text{ergodic}}$. Then, it holds that
\begin{align*}
   \hspace{-0.1cm}\frac{1}{N}\sum_{n=0}^{N-1}\mc{E}\|\nabla \bar{\mathcal{J}}(\Kbar_n)\|^2_F  &\leq \hspace{-0.05cm}\frac{2\bar{\Delta}_0}{\eta N}\hspace{-0.05cm} + \hspace{-0.05cm}\frac{c_{\text{dim}}{\epsilon}_{\text{het}}\left(\eta \hspace{-0.05cm} +\hspace{-0.05cm}\eta^2\tau_{\max}\right)}{b_s}\hspace{-0.05cm}+  \hspace{-0.05cm}c_{\text{bias}}, 
\end{align*}
where $\bar{\Delta}_0 = \mc{E}\left[\bar{\mathcal{J}}(\Kbar_{0}) - \bar{\mathcal{J}}(\Kbar^\star) \right]$, for some positive constants $c_{\text{dim}} = \mathcal{O}(n^2_x)$ and $c_{\text{bias}} = \mathcal{O}(r^2)$, and
\begin{align*}
    \eta_{\text{ergodic}} = \min \left\{\frac{1}{8}, \frac{\sqrt{b_s}}{\sqrt{6}\tau_{\max}},\frac{1}{\max\{\sqrt{32c_{\text{step}}}\tau_{\max}, 2c_{\text{step}}\}}\right\},
\end{align*}
with $c_{\text{step}} = \mathcal{O}\left(n^2_x + \frac{n^2_x}{b_s}\right)$. 
\end{theorem}

With this theorem, we are now ready to state the convergence of our Algorithm \ref{alg:async_LQR} to a first-order stationary point.
\begin{corollary} \label{cor:fixed_point_convergence_analysis} 
Let the arguments of Theorem \ref{thm:fixed_point_convergence_analysis} hold. Suppose that the step-size is set such that $\eta = \mathcal{O}\left(\sqrt{\frac{b_s}{N}}\right)$. Then, $\Kbar_N \in \mathcal{S}$ satisfy the following ergodic convergence rate\footnote{\label{footnote:bias}We omit $\mathcal{O}(c_{\text{bias}})$ in that bound, since $r$ can be set sufficiently small such that its contribution becomes negligible. See Appendix \ref{appendix:aux_lemmas}.}:
\begin{align}
   \frac{1}{N}\sum_{n=0}^{N-1}\mc{E}\|\nabla \bar{\mathcal{J}}(\Kbar_n)\|^2_F  \leq \mathcal{O}\left(\frac{\bar{\Delta}_0}{\sqrt{Nb_s}} + \frac{{\epsilon}_{\text{het}}}{\sqrt{Nb_s}} + \frac{\tau^2_{\max}\het}{N}\right). \label{eq:fixed_point}
\end{align}
\end{corollary}

Corollary~\ref{cor:fixed_point_convergence_analysis} presents the local convergence guarantee for our proposed approach with respect to the total number of iterations $N$ and batch size $b_s$. The main message of \eqref{eq:fixed_point} is that Algorithm \ref{alg:async_LQR} achieves a local convergence rate of $\mathcal{O}({\frac{1}{\sqrt{Nb_s}}})+\mathcal{O}(\frac{\tau_{\max}^2}{N}).$ In particular, the second term $\mathcal{O}(\frac{\tau_{\max}^2}{N})$ reveals the effect of the staleness, which becomes negligible when $N \geq b_s$. In the first dominant term $\mathcal{O}({\frac{1}{\sqrt{Nb_s}}})$, we demonstrate that our algorithm enjoys a linear speedup with respect to the batch size $b_s.$ This result improves upon previous work in the asynchronous distributed learning setting \cite{nguyen2022federated,toghani2022unbounded, koloskova2022sharper}, where no speedup is established.

\subsection{Global convergence guarantees}

We characterize the global convergence of the proposed asynchronous LQR design, by analyzing the interplay between staleness $\tau_{\max}$ and heterogeneity $\het$ in the optimality gap $\Delta^{(i)}_N$, i.e., the cost difference between the designed controller $\Kbar_N$ and each system's optimal controller $K^\star_i$. To this end, let us first provide an upper bound on the staleness effect throughout the iterations of Algorithm \ref{alg:async_LQR}.

\begin{lemma}\label{lemma:gradient_bound_interval} Suppose that the step-size is set according to $\eta = \mathcal{O}\left(\tau_{\max}^{-\frac{1}{2}}\right)$. Then, it holds that\footref{footnote:bias}
\begin{align*}
   \EE\|\Kbar_{l+1} - \Kbar_{l}\|^2_F \leq \eta^2 \tau_{\max}\mathcal{O}\left({\epsilon}_{\text{het}} + \EE\|\nabla \mathcal{J}^{(i)}(\Kbar_n)\|^2_F\right),
\end{align*}
$\forall l \in [n - \tau_{\max}, n-1]$ and $n \in [N-1]$.
\end{lemma}

The proof for the above lemma is detailed in Appendix \ref{appendix:convergence}. Further in this section, we present the proof sketch of our theoretical convergence guarantees, where the induction reasoning that leads to Lemma \ref{lemma:gradient_bound_interval} is discussed. By Lemma~\ref{lemma:gradient_bound_interval}, we can conclude that the staleness effect: $$\EE \lVert \Kbar_{n} -\Kstalei\rVert^2_F\le\eta^2\tau_{\max}^3\mathcal{O}\left({\epsilon}_{\text{het}} + \EE\|\nabla \mathcal{J}^{(i)}(\Kbar_n)\|^2_F \right),$$ 
can be approximately upper bounded by the product of the norm-squared policy gradient at the $n$-th iteration and the step-size squared $\eta^2$. By choosing a sufficiently small step-size $\eta$, the impact of the staleness can become negligible since it is in the high-order terms with respect to $\eta$. However, the adoption of small step-sizes $\eta$ to overcome the stragglers will slow down the convergence. We rigorously characterize such trade-off between $\tau_{\max}$ and $\eta$ in the following theorem.

\begin{theorem}[Optimality gap]\label{thm:local_optimality_gap} Let Assumptions \ref{assumption:initial_state_distribution}-\ref{assumption:bounded_het} hold. Suppose that the step-size is such that $\eta \leq \eta_{\text{gap}}$. Then, the optimality gap $\Delta^{(i)}_N = \mc{E}\left[{\mathcal{J}}^{(i)}(\Kbar_{N}) - {\mathcal{J}}^{(i)}(K^\star_i) \right]$ satisfy:
\begin{align*}
     \Delta^{(i)}_N \leq c_{\text{cont}}^N\Delta^{(i)}_0  +  \lambda^{-1}(c_{\text{dim}}{\epsilon}_{\text{het}} + c_{\text{bias}}),
\end{align*}
where $c_{\text{cont}} = 1 - \frac{\eta\lambda}{4} \in (0,1)$ denotes the contraction rate and $ \eta_{\text{gap}} = \min \left\{\eta_{\text{ergodic}}, \frac{1}{\barhgrad \max\left\{8\tau_{\max},\sqrt{2}\tau^{3/2}_{\max}\right\}},\frac{1}{\tau_{\max}\sqrt{2c_{\text{bias}}}}\right\}$.
\end{theorem}

As stated in the above theorem, the optimality gap $\Delta^{(i)}_N$, is composed of a contraction term, with contraction rate $c_{\text{cont}}$, and an additive bias characterized by $\het$ and $c_{\text{bias}}$. As the number of iterations $N$ grows, the contraction shrinks to zero. In addition, as $c_{\text{bias}}$ is in the order of $r^2$, the smoothing radius can be set sufficiently small such that its contribution becomes negligible in the bias term. However, $\het$ is fixed and dominates the unavoidable bias in the optimality gap. In addition, the step-size $\eta$ is in the order of $\mathcal{O}\left(\frac{1}{\tau^{3/2}_{\max}}\right)$. This condition demonstrates that as the staleness $\tau_{\max}$ increases, the step-size $\eta$ needs to be reduced to preserve a global convergence guarantee of Algorithm \ref{alg:async_LQR}. Next, we highlight the impact of $\tau_{\max}$ on the number of iterations $N$ to achieve a controller that is $\epsilon$-close to its optimal controller. 

\begin{corollary}[Linear convergence]\label{corollary:linear_convergence} Let the arguments of Theorem \ref{thm:local_optimality_gap} hold. Suppose that the number of iterations $N$ of Algorithm \ref{alg:async_LQR} and the smoothing radius $r$ of Algorithm \ref{alg:ZO2P} satisfy:
$N \geq \mathcal{O}\left(\tau^{3/2}_{\max}\log\left(\frac{\Delta^\vee_0}{\epsilon}\right)\right), r \leq \mathcal{O}(\epsilon),$ for a small tolerance $\epsilon \in (0,1)$, where $\Delta^\vee_0 := \max_{i \in [M]} \Delta^{(i)}_0$. Then, the optimality gap satisfies: $\Delta^{(i)}_N \leq \mathcal{O}\left( \epsilon  + {\epsilon}_{\text{het}}\right)$.
\end{corollary}

Corollary \ref{corollary:linear_convergence} shows that by carefully controlling the step-size $\eta$, number of iterations $N$ and smoothing radius $r$, the designed controller $\Kbar_N$ is $\epsilon$-close to each system's optimal controllers up to a heterogeneity bias. Note that the number of iterations $N$ will increase with the maximum staleness -- which is of order  $\tau_{max}^{3/2}.$ In other words, the staleness will slow down the global convergence rate.

\subsection{Proof sketch}

We now discuss the main steps and reasoning to obtain the theoretical convergence results presented in this work. First, for the local convergence rate, we note that as long as Lemma \ref{lemma:smoothness} holds for any system $i \in [M]$, it implies that the gradient of the average LQR cost is also $\barhgrad$-Lipschitz. Therefore, the gap in the average cost between two consecutive iterations of Algorithm \ref{alg:async_LQR}, i.e.,  $\bar{\Delta}_n = \bar{\mathcal{J}}(\Kbar_{n+1}) - \bar{\mathcal{J}}(\Kbar_{n})$, is approximately upper bounded as follows 
\begin{align*}
    \mc{E}\left[\bar{\Delta}_n \right] &\lesssim -\eta \mc{E}\|\nabla \bar{\mathcal{J}}(\Kbar_n)\|^2_F - \eta \mc{E}\left\| \nabla \bar{\mathcal{J}}_r(\Kstalei)\right\|_F^2\\ &+\frac{\eta}{M}\sum_{i=1}^M  \underbrace{\mc{E} \left\| \Kbar_n - \Kstalei\right\|_F^2}_{\text{staleness term}}
    + \eta( r^2 + \eta\het), 
\end{align*}
where we use $\lesssim$ to omit constant factors in the expression, and $\nabla\bar{\mathcal{J}}_r(\Kbar) := \mc{E}\widehat{\nabla}\bar{\mathcal{J}}(\Kbar)$. In addition, the staleness term can be upper bounded as follows
\begin{align}\label{eq:staleness_term}
\mc{E} \left\| \Kbar_n - \Kstalei\right\|_F^2 &= \mc{E} \left\|\sum_{l = n - \tau_i(n)}^{n-1} \Kbar_{l+1} - \Kbar_{l}\right\|_F^2\notag\\
&\leq  \tau_{\max}\sum_{l = n - \tau_i(n)}^{n-1}\mc{E} \left\| \Kbar_{l+1} - \Kbar_{l}\right\|_F^2,
\end{align}
then, by summing the above expression over the iterations $n$, the staleness effect is shown to be in the following order.
\begin{align*}
    \hspace{-0.2cm}\sum_{n = 0}^{N-1}\mc{E} \left\| \Kbar_n - \Kstalei\right\|_F^2 &\lesssim \tau^2_{\max}N\eta^3r^2 + \frac{\tau^2_{\max}N\eta^3\het}{b_s}\\
    &\hspace{-0.4cm}+\tau^2_{\max}\eta^3\sum_{n=0}^{N-1}\mc{E} \left\| \nabla \bar{\mathcal{J}}_r (\Kstalei)\right\|^2_F,
\end{align*}
which can be used in the expression of the expected average gap, i.e.,  $\mc{E}\left[\bar{\Delta}_n \right]$, and with a proper selection of the step-size $\eta$ we can obtain the result presented in Theorem \ref{thm:fixed_point_convergence_analysis}. We emphasize that since the interplay between staleness $\tau^3_{\max}$ and heterogeneity $\het$ is accompanied by $\eta^3$, its contribution in the local convergence rate should only appears in a high-order term when we set $\eta = \mathcal{O}\left(\sqrt{\frac{b_s}{N}}\right)$. 

On the other hand, what changes in the global convergence guarantees is how the staleness effect is upper bounded. To this end, we use the result in Lemma \ref{lemma:gradient_bound}. The proof idea of this lemma relies on an induction approach, where for any two consecutive iterations of Algorithm \ref{alg:async_LQR}, we can write  
\begin{align}\label{eq:staleness_all_n}
    \mc{E}\|\Kbar_{n+1} - \Kbar_{n}\|^2_F \lesssim \tau_{\max}\eta^2\het  + \tau_{\max}\eta^2\mc{E}\|\nabla \mathcal{J}^{(i)}(\Kbar_n)\|^2_F,
\end{align}
then, since the staleness term \eqref{eq:staleness_term} is evaluated within the interval $l \in [n - \tau_{\max}, n-1]$, the bound in \eqref{eq:staleness_all_n} combined with an induction step implies in Lemma \ref{lemma:gradient_bound_interval}. Therefore, by using such lemma with the local smoothness and gradient dominance properties (i.e., Lemmas \ref{lemma:smoothness}-\ref{lemma:gradient_dominance}), we obtain the global convergence results of Theorem \ref{thm:local_optimality_gap} and Corollary \ref{corollary:linear_convergence}.

\subsection{Stability guarantees}

An important requirement that policy gradient methods needs to satisfy within control tasks, is the ability of iteratively preserving the closed-loop stability of the collaborating systems with respect to the designed controller. Note that, one of the conditions to ensure such requirement, is to have access to an initial stabilizing controller $\Kbar_0 \in \mathcal{S}$. However, it is also necessary to impose conditions on the step-size $\eta$, smoothing radius $r$, and heterogeneity $\het$  to ensure that big steps, non-accurate PG estimates, and large heterogeneous settings\footnote{The authors in \cite{wang2023model} discuss the necessity of a low heterogeneity regime when designing stabilizing LQR controllers in a multi-agent setting.} will not produce unstabilizing policy updates. We summarize such conditions in the following theorem.

\begin{theorem}[Per-iteration stabilizing controllers]\label{thm:stability_analysis} Let Assumptions \ref{assumption:initial_state_distribution}-\ref{assumption:bounded_het} hold. Suppose that the step-size is set such that $\eta \leq\eta_{\text{gap}}$. In addition, suppose that the heterogeneity and smoothing radius satisfy $\het \leq \frac{\gamma \lambda \Delta^\vee_0}{64}$ and $r^2 \leq\frac{\gamma \lambda \Delta^\vee_0}{64\barhgrad^2}$, respectively. Then, Algorithm \ref{alg:async_LQR} produces a stabilizing controller $\Kbar_n \in \mathcal{S}$ for all iterations $n \in \{0,1,\ldots,N-1\}$.
\end{theorem}

The proof for this theorem is detailed in Appendix \ref{appendix:stability}, where given an initial stabilizing controller, we exploit an induction approach along with Definition \ref{def:stabilizing_set} to derive the necessary conditions on the step-size $\eta$, smoothing radius $r$ and heterogeneity $\het$ that ensure the design of stabilizing controllers $\Kbar_n$, for all iterations $n$ of Algorithm \ref{alg:async_LQR}.  

\section{Numerical Experiments} \label{sec:numerical_results}

Numerical experiments\footnote{Code to reproduce the results can be downloaded from: \url{https://github.com/jd-anderson/AsyncLQR}.} are now presented to illustrate and validate the convergence guarantees of Algorithm \ref{alg:async_LQR}. In particular, we highlight the effect of the staleness $\tau_{\max}$, batch size $b_s$ and heterogeneity level $\het$ to the speed of convergence and optimality gap of Algorithm \ref{alg:async_LQR}. We also illustrate the benefit of asynchronous aggregation over the synchronous counterpart \cite{wang2023model}, when dealing with straggler systems in the learning process. To this end, let us first consider \emph{nominal} system matrices: 
\begin{align*}
\underbrace{\left[\begin{array}{cccc}1.22 & 0.03 & -0.02 & -0.32 \\
0.01 & 0.47 & 4.70 & 0.00 \\
0.02 & -0.06 & 0.40 & 0.00 \\
0.01 & -0.04 & 0.72 & 1.55
\end{array}\right]}_{A^{(1)}}, 
\underbrace{
\left[\begin{array}{cc}0.01 & 0.99 \\
-3.44 & 1.66 \\
-0.83 & 0.44 \\
-0.47 & 0.25
\end{array}\right]}_{B^{(1)}},
\end{align*}
and cost matrices $Q^{(1)} = I_{4}$ and $R^{(1)} = I_2$. Therefore, by applying random perturbations to $(A^{(1)}, B^{(1)}, Q^{(1)}, R^{(1)})$, with radius $(\epsilon_A, \epsilon_B, \epsilon_Q, \epsilon_R)$, we generate $M = 100$ tuples $(A^{(i)}, B^{(i)}, Q^{(i)}, R^{(i)})$, for $i \in [M]$, to characterize our heterogeneous multi-agent LQR setting. We provide more details on that generative process in Appendix \ref{appendix:numericals}.

\begin{figure}
    \centering
    \includegraphics[width=0.5\textwidth]{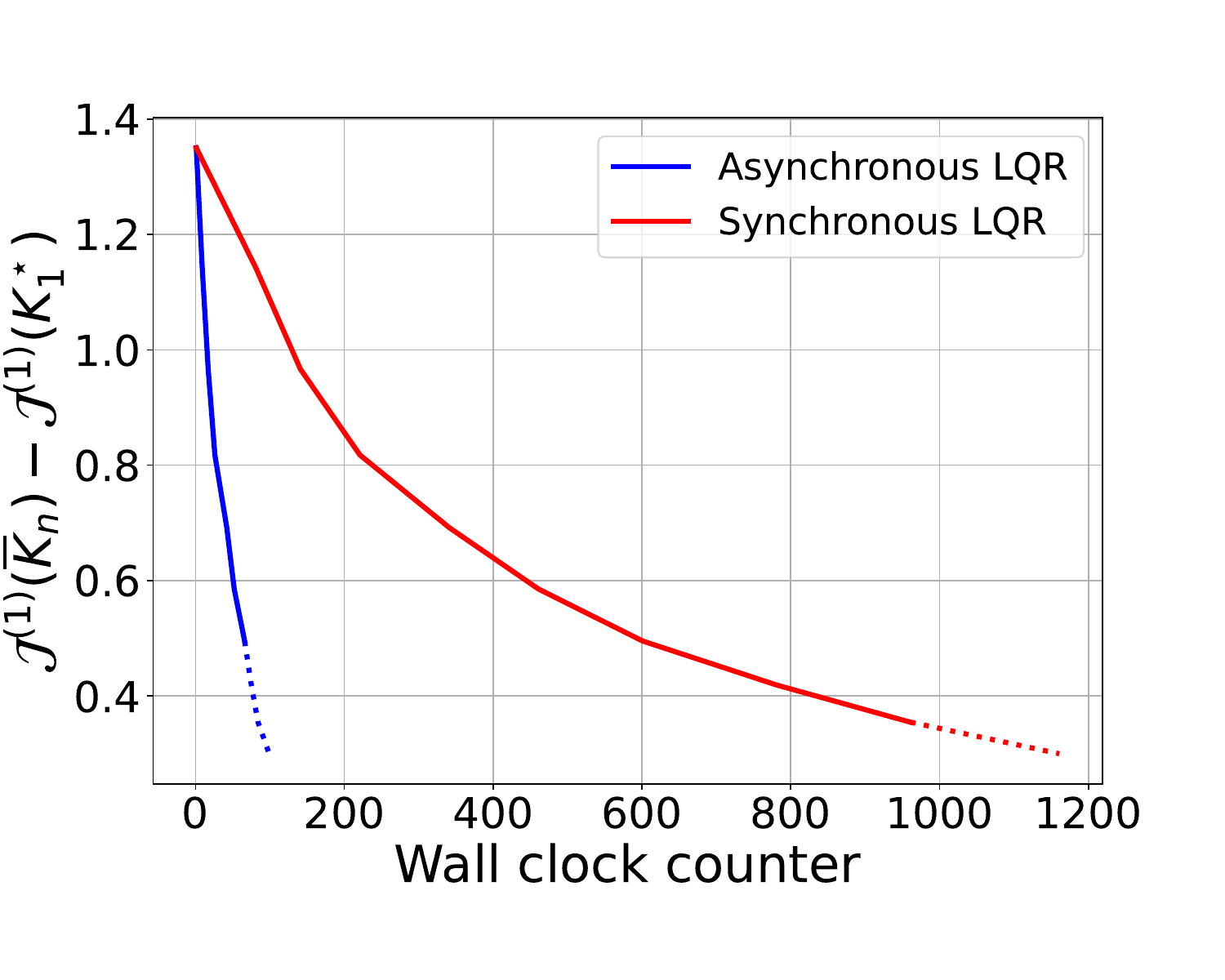}\vspace{-0.8cm}
    \caption{Optimality gap (with respect to the nominal system and cost matrices) as a function of the wall clock counter, for the asynchronous LQR design (Algorithm \ref{alg:async_LQR}) and the synchronous LQR approach \cite{wang2023model}, both in the presence of a single straggler system with $\tau_{\max} = 20$. We set $(\epsilon_A = 5.46,  \epsilon_B = 2.74, \epsilon_Q = 3.96 , \epsilon_R = 2.82)\times 10^{-2}$ and  $b_s = 20$.} 
    \label{fig:comparison_sync_async}
\end{figure}

\begin{figure*}
    \centering
   \subfigure[Varying $\tau_{\max}$]{\includegraphics[width=0.33\textwidth]{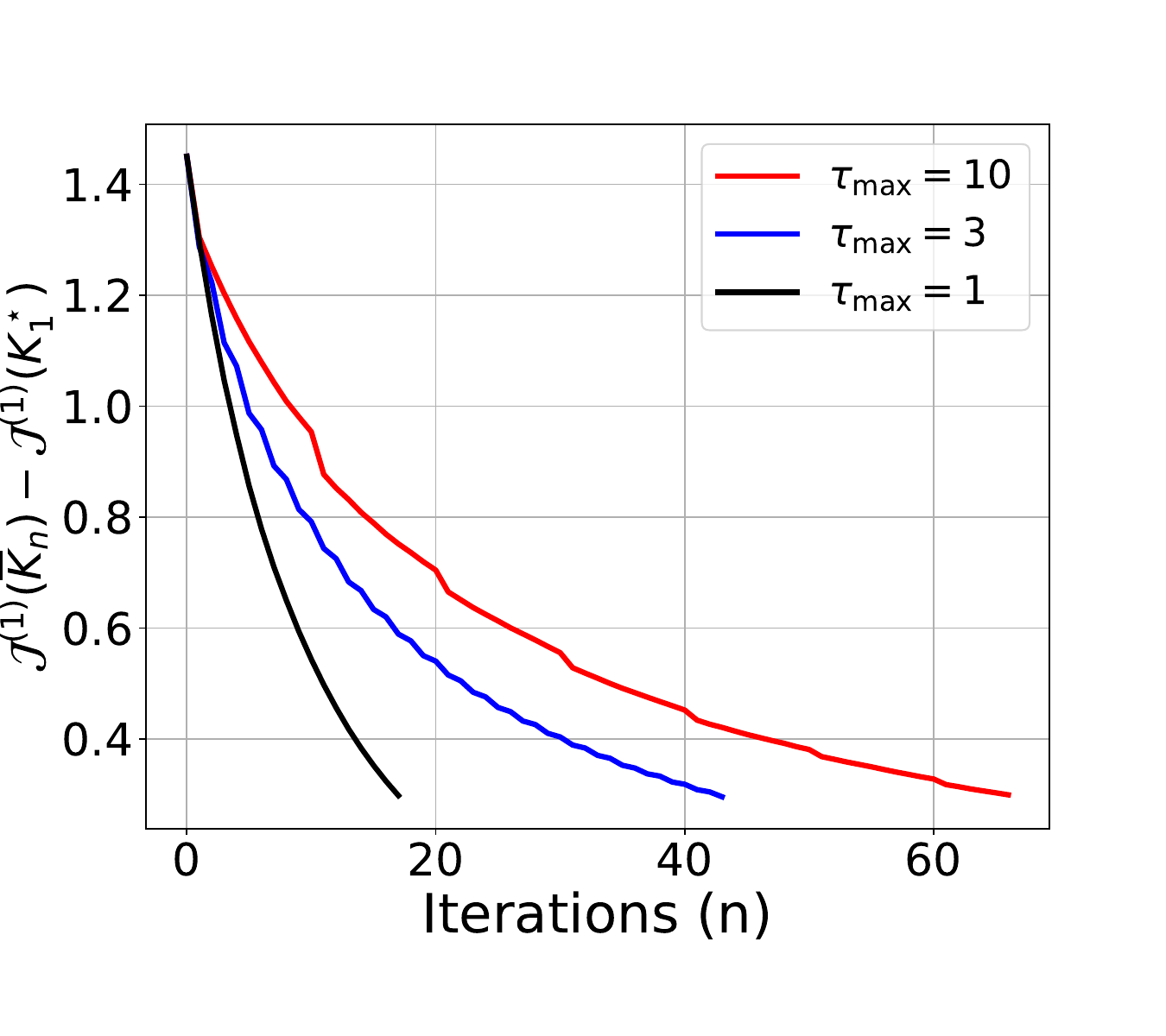}}\label{fig3} \hspace{-0.6cm}
   \subfigure[Varying $b_s$]{\includegraphics[width=0.33\textwidth]{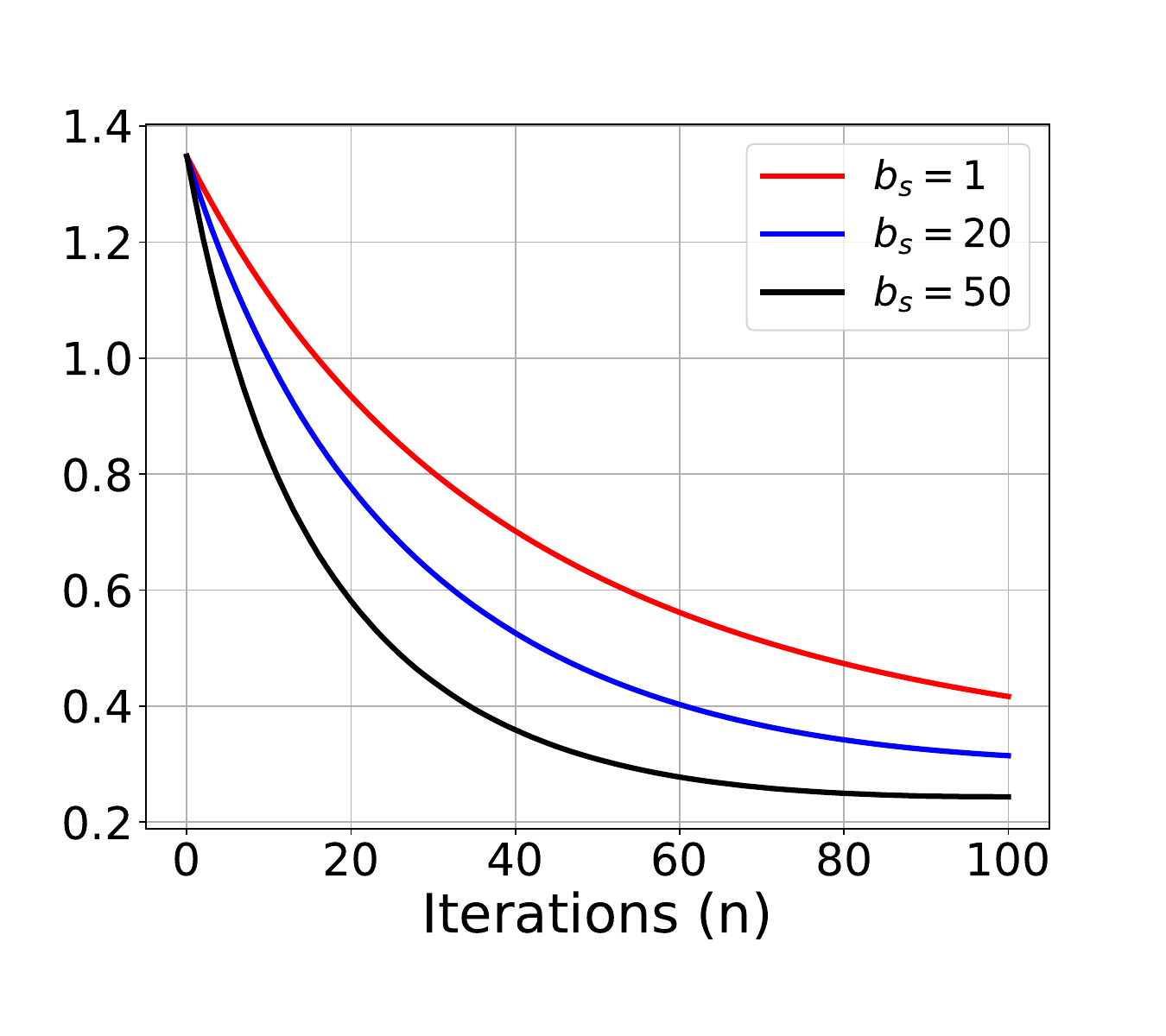}}\label{fig2}\hspace{-0.4cm}
   \subfigure[Varying $\het$]{\includegraphics[width=0.33\textwidth]{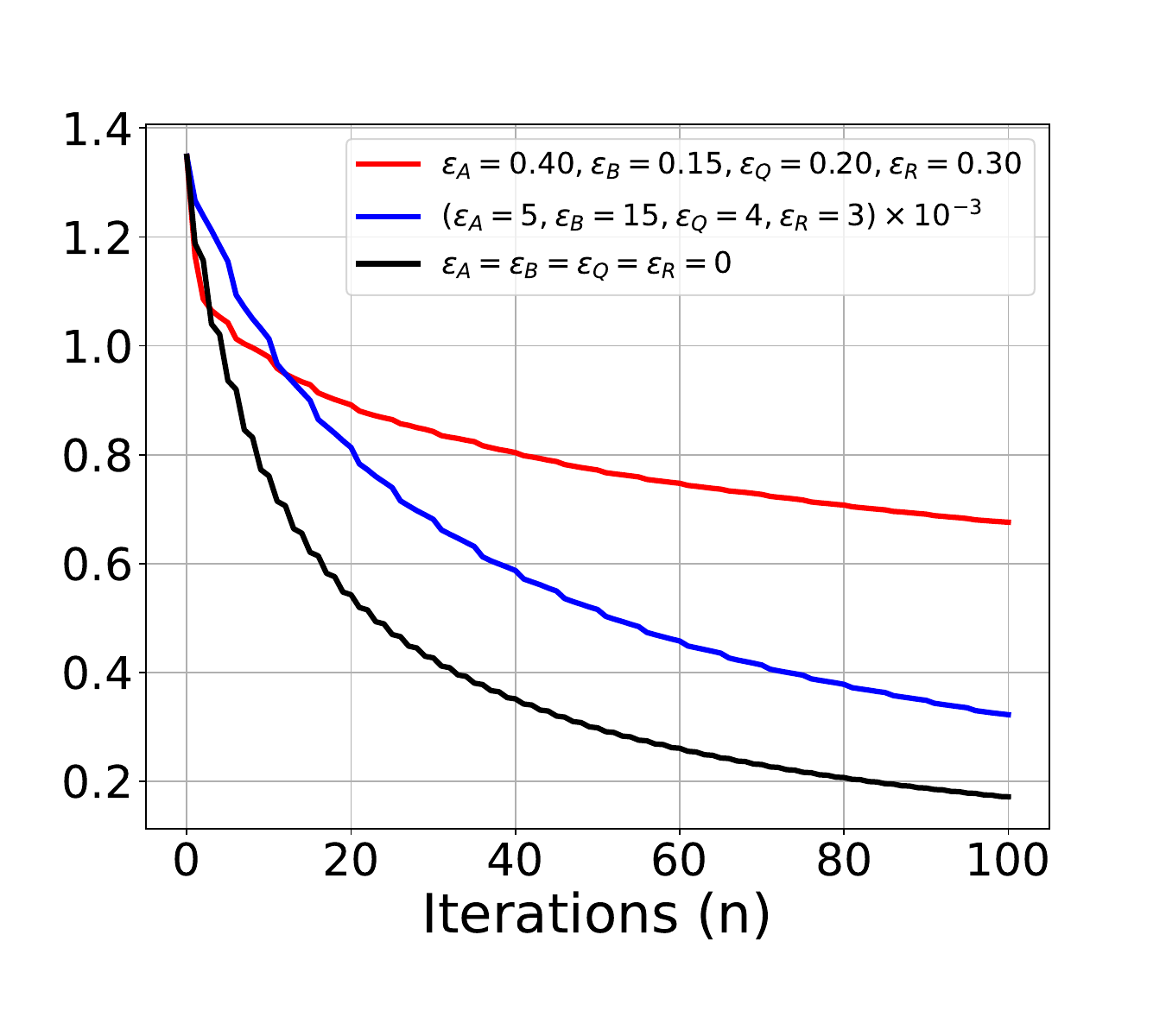}}\label{fig4}\hspace{-0.5cm}
   \vspace{-0.2em}
    \caption{Optimality gap with respect to the iteration count $n$. (a) $b_s = 20$ and $(\epsilon_A = 5.46,  \epsilon_B = 2.74, \epsilon_Q = 3.96 , \epsilon_R = 2.82)\times 10^{-2}$. (b) $\tau_{\max} = 1$ and $(\epsilon_A = 5.25, \epsilon_B = 2.80, \epsilon_Q = 4.00, \epsilon_R = 2.82)\times 10^{-3}$. (c) $\tau_{\max} = 5$ and $b_s = 20$.}
    \label{fig:convergence}
\hrulefill
\end{figure*}

With $M$ system and cost matrices in hands, we first compare the proposed asynchronous LQR design of Algorithm \ref{alg:async_LQR} over the synchronous federated LQR approach in \cite{wang2023model}, both in the presence of a single straggler system with $\tau_{\max} = 20$.  Since, in Algorithm \ref{alg:async_LQR}, the server performs controller updates upon receiving the fastest $b_s$ PG estimates, such quicker systems are not affected by straggler systems when $\tau_{\max}$ is sufficiently large. To illustrate this, Figure \ref{fig:comparison_sync_async} shows how long, in terms of a wall clock counter, Algorithm \ref{alg:async_LQR} takes to design a controller $\Kbar_N$ that achieves a certain optimality gap compared to the synchronous LQR design. This figure shows that due to the presence of stragglers, the synchronous federated LQR approach \cite{wang2023model} needs to \emph{wait} a long time for all of the $M$ PG estimates to be reported to the server to then proceed with the controller update. On the other hand, Algorithm \ref{alg:async_LQR} fully enjoys the parallelism of distributed computation, even when dealing with slow systems. However, as discussed, this comes with the price of aggregating stale PG estimates. Figure \ref{fig:convergence} illustrates the effect of $\tau_{\max}$ in the optimality gap of Algorithm \ref{alg:async_LQR}.

%One paragraph discussing the ergodic rate results
 
Figure \ref{fig:convergence} depict the optimality gap $\Delta^{(1)}_n$ as a function of the iteration count $n$, for a varying: (a) staleness $\tau_{\max}$, (b) batch size $b_s$, and (c) heterogeneity level $\het$. Note that, for convenience, we evaluate the global convergence of Algorithm \ref{alg:async_LQR} on the nominal system and cost, i.e., $i=1$. However, we emphasize that a  similar result should also be observed for any $i \in [M]$. Figure \ref{fig:convergence}-(a) shows the impact of $\tau_{\max}$ in the number of iterations $N$ needed to achieve a certain optimality gap; it highlights that the number of iterations required to design a controller $\Kbar_N$ such that $\Delta^{(1)}_N \leq 0.3$ is larger when $\tau_{\max} = 3$ and $\tau_{\max} = 10$ compared to $\tau_{\max} = 1$. Moreover, as predicted in Corollary \ref{corollary:linear_convergence}, the staleness $\tau_{\max}$ only affects the speed of convergence and does not impact the accuracy in the optimality gap $\Delta^{(i)}_N$.

Furthermore, in alignment with Corollary \ref{cor:fixed_point_convergence_analysis}, Figure \ref{fig:convergence}-(b) illustrates the benefit of aggregating multiple system's PG estimates. As predicted, an increase in the batch size $b_s$ leads to a faster convergence of Algorithm \ref{alg:async_LQR}. Lastly, as illustrated in Figure \ref{fig:convergence}-(c), due to the heterogeneous setting, Algorithm \ref{alg:async_LQR} returns an $\epsilon$-near optimal controller up to a heterogeneity bias. Therefore, as $\het$ increases, the unavoidable bias (Corollary \ref{corollary:linear_convergence}) also increases.

\section{Conclusions and Future Work} \label{sec:conclusions}

To understand how aggregating stale policy gradient estimates affect model-free LQR design, we characterized the convergence and stability guarantees of an asynchronous and heterogeneous PG method applied to the multi-agent LQR problem. Despite straggler systems, the proposed asynchronous aggregation scheme fully exploits the parallelism in the distributed computation (see Figure \ref{fig:comparison_sync_async}). Nevertheless, such parallelism comes with the price of aggregating stale policy gradient estimates. Our analysis demonstrated that, by carefully controlling the step-size, the staleness effect remains limited to a high-order term of the ergodic convergence rate (Corollary \ref{cor:fixed_point_convergence_analysis}). Moreover, the optimality gap bound remains untouched as in the synchronous case \cite{wang2023model}. We showed that the staleness impacts the speed of convergence through a multiplicative factor (Corollary \ref{corollary:linear_convergence}). We provided numerical results to illustrate and validate our theory (i.e., Figure \ref{fig:convergence}), where we also highlight the effect of the heterogeneity level $\het$ and batch size $b_s$ to the convergence of Algorithm \ref{alg:async_LQR}. Future work may explore other aggregation schemes, beyond a simple average, to alleviate the staleness effect even more in the local and global convergence bounds. 

\bibliography{references.bib}
\bibliographystyle{IEEEtran}
\newpage

\onecolumn

\appendix

In this appendix, we provide the technical details and proofs for the convergence and stability guarantees of Algorithm \ref{alg:async_LQR}. For this purpose, we first introduce some auxiliary lemmas that are instrumental in deriving the results of Theorems \ref{thm:fixed_point_convergence_analysis}-\ref{thm:local_optimality_gap} and Corollaries \ref{cor:fixed_point_convergence_analysis}-\ref{corollary:linear_convergence}. In particular, Lemma \ref{lemma:gradient_bound} highlights that with a careful selection of the step-size $\eta  = \mathcal{O}(\tau^{-1/2}_{\max})$, an induction step suffices to show that even in the presence of staleness, $\mc{E}\|\Kbar_{n+1} - \Kbar_{n}\|^2_F$ is in the order of the $\mc{E}\|\nabla \mathcal{J}^{(i)}(\Kbar_n)\|^2_F$ up to a heterogeneity and bias term. With a proper analysis, that result implies in the staleness upper bound in Lemma \ref{lemma:gradient_bound_interval}. This result is crucial to show the linear convergence in the optimality gap (i.e., Theorem \ref{thm:local_optimality_gap} and Corollary \ref{corollary:linear_convergence}). Sections \ref{appendix:stability} and \ref{appendix:convergence} details the local and global convergence, and stability analysis of Algorithm \ref{alg:async_LQR}, respectively. Moreover, in Section \ref{appendix:numericals}, we provide additional details on the experimental setup we used to produce Figures \ref{fig:comparison_sync_async}-\ref{fig:convergence}.

\subsection{Auxiliary Lemmas}\label{appendix:aux_lemmas}

\begin{lemma}[Young's inequality]\label{lemma:youngs}
Given two matrices $A, B \in \mathbb{R}^{n_x\times n_u}$,  for any $\alpha>0$, we have

\begin{align*}%\label{eq:youngs}
\|A+B\|_2^2 \leq (1+\alpha)\|A\|_F^2+\left(1+\frac{1}{\alpha}\right)\|B\|_F^2.
\end{align*}

Furthermore, 

\begin{align*}
\langle A, B\rangle  \leq  \frac{\alpha}{2}\lVert A\rVert_F^2 +\frac{1}{2\alpha}\lVert B \rVert_F^2.
\end{align*}
\end{lemma}

\begin{lemma}[Jensen's inequality]\label{lemma:jensen}
 Given a collection of $M$ matrices $\{A^{(i)}\}_{i=1}^{M}$ with identical dimensions, it holds that 
\begin{align*}
\left\|\sum_{i=1}^M A^{(i)}\right\|_F^2 \leq M \sum_{i=1}^M\left\|A^{(i)}\right\|_F^2.
\end{align*}
\end{lemma}

\begin{lemma}[Controlling the bias, Lemma 5 of \cite{malik2019derivative}] \label{lemma:bias}  Let $\widehat{\nabla}\mathcal{J}(K)$ be the two-point zeroth-order gradient estimations evaluated at a stabilizing controller $K \in \mathcal{S}$. Then, it holds that
\begin{align*}%\label{eq:control_bias}
    \mc{E}\|\nabla \mathcal{J}(K) - {\nabla}J_r(K)\|^2_F \leq  h^2_{\text{grad}} r^2.
\end{align*}
where ${\nabla}J_r(K) \triangleq \mc{E}\widehat{\nabla}\mathcal{J}(K)$.
\end{lemma}

\begin{lemma}[Lemma 4 of \cite{toso2023oracle}] \label{lemma:two_point_estimated_grad}
Let $\widehat{\nabla}\mathcal{J}(K)$ be the two-point zeroth-order gradient estimation of ${\nabla}\mathcal{J}(K)$ evaluated at a stabilizing controller $K$. Then, it holds that 
\begin{align*}
    \mc{E}\left\|\widehat{\nabla}\mathcal{J}(K)\right\|^2_F  \leq c_{\text{ZO}}(n_x)\left(
    \barhgrad^2r^2+\mc{E}\|\nabla \mathcal{J}(K) \|^2_F \right),
\end{align*}
with $c_{\text{ZO}}(n_x) = 8n_x^2$.
\end{lemma}

\begin{lemma}[Lemma 7 of \cite{cheng2024momentum}]\label{lem:par_sample}
Given vectors $\{a_1, \cdots, a_N, b\}\in \mathbb{R}^d$, of dimension $d$ with $a=\frac{1}{N}\sum_{i\in[N]}a_i$ and uniform sampling $\mathcal{G}\subset [N]$ without replacement such that $|\mathcal{G}|=G$, it holds that
    \begin{align*}%\label{eq:par_sample}
            \mc{E}_{\mathcal{G}} \left[\left\|\frac{1}{G}\sum_{i\in \mathcal{G}}a_i\right\|^2 \right]
            \leq  &\|  a\|^2 + \frac{1}{GN}\sum_{i}\|a_i-{a}\|^2\leq \|  a\|^2 + \frac{1}{GN}\sum_{i}\|a_i-b\|^2.         
 \end{align*}
\end{lemma}

\begin{lemma}\label{lemma:gradient_bound} Suppose that the step-size is selected according to $\eta \leq \frac{c}{\barhgrad \sqrt{\max\{c^1_{\tau_{\max}},c^2_{\tau_{\max}},\tau_{\max}\}}}$, for a sufficiently small constant $c \in \mathbb{R}_{>0}$, with $c^1_{\tau_{\max}}, c^2_{\tau_{\max}} = \mathcal{O}(\tau_{\max})$. Then, for any iteration $n \in \{0,1,\ldots, N-1\}$ of Algorithm \ref{alg:async_LQR}, it holds that
\begin{align}\label{eq:bound_on_the_gradient}
    \mc{E}\|\Kbar_{n+1} - \Kbar_{n}\|^2_F \leq c^1_{\tau_{\max}}\eta^2{\epsilon}_{\text{het}}  + c^2_{\tau_{\max}}\eta^2\mc{E}\|\nabla \mathcal{J}^{(i)}(\Kbar_n)\|^2_F + \eta^2c_{\text{bias}},
\end{align}
for some positive constant $c_{\text{bias}} = \mathcal{O}(r^2)$.

\begin{proof}
We begin our proof by using the updating rule in Algorithm \ref{alg:async_LQR} to write  
\allowdisplaybreaks
\begin{align*}
   \mc{E} \|\Kbar_{n+1} &- \Kbar_n\|^2_F \leq  \frac{\eta^2}{b_s}\sum_{s=1}^{b_s}\mc{E}\left\|\widehat{\nabla} \mathcal{J}^{(s)}(\Kbar_{n-\tau_s(n)})\right\|^2_F\notag\\
    &\stackrel{(i)}{\leq}  \frac{\eta^2}{b_s}\sum_{s=1}^{b_s}\mc{E}\left\|{\nabla} \mathcal{J}^{(s)}(\Kbar_{n-\tau_s(n)})\right\|^2_F + \eta^2c_{\text{ZO}}(n_x)\barhgrad^2 r^2\notag\\
    &= \eta^2 \frac{1}{b_s}\sum_{s=1}^{b_s}\mc{E}\left\|\nabla \mathcal{J}^{(s)}(\Kbar_{n-\tau_s(n)}) - \nabla \mathcal{J}^{(i)}(\Kbar_n) + \nabla \mathcal{J}^{(i)}(\Kbar_n) \right\|^2_F + \eta^2c_{\text{ZO}}(n_x)\barhgrad^2 r^2\notag\\
    &\stackrel{(ii)}{\leq} 4\eta^2 \frac{1}{b_s}\sum_{s=1}^{b_s}\mc{E}\left\|\nabla \mathcal{J}^{(s)}(\Kbar_{n-\tau_s(n)}) - \nabla \mathcal{J}^{(s)}(\Kbar_n)\right\|^2_F +4\eta^2 {\epsilon}_{\text{het}} + 2\eta^2 \mc{E}\left\|\nabla \mathcal{J}^{(i)}(\Kbar_n)\right\|^2_F + \eta^2c_{\text{ZO}}(n_x)\barhgrad^2 r^2\notag\\
    &\stackrel{(iii)}{\leq} 4\eta^2\barhgrad^2 \frac{1}{b_s}\sum_{s=1}^{b_s}\mc{E}\left\| \Kbar_{n-\tau_s(n)} - \Kbar_n\right\|^2_F +4\eta^2 {\epsilon}_{\text{het}} + 2\eta^2 \mc{E}\left\|\nabla \mathcal{J}^{(i)}(\Kbar_n)\right\|^2_F + \eta^2c_{\text{ZO}}(n_x)\barhgrad^2 r^2\notag\\
    &\leq 4\tau_{\max}\eta^2\barhgrad^2 \sum_{l=n-\tau_{\max}}^{n-1} \mc{E}\left\|\Kbar_{l+1} - \Kbar_{l}\right\|^2_F +4\eta^2 {\epsilon}_{\text{het}} + 2\eta^2 \mc{E}\left\|\nabla \mathcal{J}^{(i)}(\Kbar_n)\right\|^2_F + \eta^2c_{\text{ZO}}(n_x)\barhgrad^2 r^2,
\end{align*}
for any iteration $n \in \{0,1,\ldots,N-1\}$. In $(i)$ we use Lemma \ref{lemma:two_point_estimated_grad}, whereas $(ii)$ follows from the gradient heterogeneity bound in Lemma \ref{lemma:gradient_heterogeneity}. In $(iii)$ we use the Lipschitz property of the LQR gradient. We then proceed to show \eqref{eq:bound_on_the_gradient}, by using an induction approach. We first define the base case as the first iteration $n=0$ which yields 
\begin{align*}
   \textbf{Base case:}\;\ \mc{E}\|\Kbar_{1} - \Kbar_0\|^2_F &\leq  4\eta^2 {\epsilon}_{\text{het}} + 2\eta^2 \mc{E}\left\|\nabla \mathcal{J}^{(i)}(\Kbar_0)\right\|^2_F + c_{\text{ZO}}(n_x)\eta^2\barhgrad^2 r^2\\
   &\leq  c^1_{\tau_{\max}}\eta^2 {\epsilon}_{\text{het}} + c^2_{\tau_{\max}}\eta^2 \mc{E}\left\|\nabla \mathcal{J}^{(i)}(\Kbar_0)\right\|^2_F + \eta^2c_{\text{bias}},
\end{align*}
for some positive constants $c^1_{\tau_{\max}}$, $c^2_{\tau_{\max}}$ in the order of $\tau_{\max}$. Then, our inductive hypothesis is defined as follows
\begin{align*}
   \textbf{Inductive hypothesis:}\;\ \mc{E}\|\Kbar_{n} - \Kbar_{n-1}\|^2_F \leq c^1_{\tau_{\max}}\eta^2 {\epsilon}_{\text{het}} + c^2_{\tau_{\max}}\eta^2 \mc{E}\left\|\nabla \mathcal{J}^{(i)}(\Kbar_{n-1})\right\|^2_F + \eta^2c_{\text{bias}}, \text{ $\forall n \in [N-1]$.}
\end{align*}

Note that for any iteration $n+1 \in \{0,1,\ldots,N-1\}$ we have 
\begin{align}\label{eq:update_rule_gradient_bound}
    \mc{E}\|\Kbar_{n+1} - \Kbar_n\|^2_F \leq 4\tau_{\max}\eta^2\barhgrad^2 \sum_{l=n-\tau_{\max}}^{n-1} \mc{E}\left\|\Kbar_{l+1} - \Kbar_{l}\right\|^2_F +4\eta^2 {\epsilon}_{\text{het}} + 2\eta^2 \mc{E}\left\|\nabla \mathcal{J}^{(i)}(\Kbar_n)\right\|^2_F +  c_{\text{ZO}}(n_x)\eta^2\barhgrad^2 r^2.
\end{align}

We then proceed by showing the bound in \eqref{eq:bound_on_the_gradient} holds for any $\tau_{\max} \leq n$. First, for $\tau_{\max} = 1$, we have 

\begin{align}\label{eq:gradient_bound_tau_1}
    \mc{E}\|\Kbar_{n+1} - \Kbar_n\|^2_F &\leq 4\tau_{\max}\eta^2\barhgrad^2 \mc{E}\left\|\Kbar_{n} - \Kbar_{n-1}\right\|^2_F + 4\eta^2 {\epsilon}_{\text{het}} + 2\eta^2 \mc{E}\left\|\nabla \mathcal{J}^{(i)}(\Kbar_n)\right\|^2_F +  c_{\text{ZO}}(n_x)\eta^2\barhgrad^2 r^2\notag\\
    &\stackrel{(i)}{\leq} 4c^1_{\tau_{\max}}\tau_{\max}\eta^4\barhgrad^2{\epsilon}_{\text{het}} + 4\eta^2 {\epsilon}_{\text{het}} + 4c^2_{\tau_{\max}}\tau_{\max}\eta^4\barhgrad^2\mc{E}\left\|\nabla \mathcal{J}^{(i)}(\Kbar_{n-1})\right\|^2_F + 2\eta^2 \mc{E}\left\|\nabla \mathcal{J}^{(i)}(\Kbar_n)\right\|^2_F\notag\\
    &+4\tau_{\max}\eta^4\barhgrad^2c_{\text{ZO}}(n_x)\barhgrad^2r^2 +  c_{\text{ZO}}(n_x)\eta^2\barhgrad^2 r^2\notag\\
    &\stackrel{(ii)}{\leq} 4(1+\tau_{\max})\eta^2{\epsilon}_{\text{het}} +8\tau_{\max}\eta^2\barhgrad^2\|\Kbar_n-\Kbar_{n-1}\|^2_F + 2(1+4\tau_{\max})\eta^2 \mc{E}\left\|\nabla \mathcal{J}^{(i)}(\Kbar_n)\right\|^2_F\notag\\
    &+2c_{\text{ZO}}(n_x)\eta^2\barhgrad^2 r^2,
\end{align}
where $(i)$ follows from the inductive hypothesis, whereas in $(ii)$ follows from the  selection of the the step-size. We then note that for $n = 1$, we obtain
\begin{align}\label{eq:gradient_bound_tau_1_n_1}
    \mc{E}\|\Kbar_{2} - \Kbar_1\|^2_F &\leq 4(1+\tau_{\max})\eta^2{\epsilon}_{\text{het}} +8\tau_{\max}\eta^2\barhgrad^2\mc{E}\|\Kbar_{1} - \Kbar_0\|^2_F + 2(1+4\tau_{\max})\eta^2 \mc{E}\left\|\nabla \mathcal{J}^{(i)}(\Kbar_1)\right\|^2_F\notag\\
    &+2c_{\text{ZO}}(n_x)\eta^2\barhgrad^2 r^2\notag \\
    &\stackrel{(i)}{\leq} c^1_{\tau_{\max}}\eta^2{\epsilon}_{\text{het}} +  c^2_{\tau_{\max}}\eta^2 \mc{E}\left\|\nabla \mathcal{J}^{(i)}(\Kbar_1)\right\|^2_F +\eta^2c_{\text{bias}},
\end{align}
where $(i)$ follows from the base case and the condition on the step-size. By setting $n=2$ in \eqref{eq:gradient_bound_tau_1}, we have
\begin{align*}
    \mc{E}\|\Kbar_{3} - \Kbar_2\|^2_F &\leq 4(1+\tau_{\max})\eta^2{\epsilon}_{\text{het}} +8\tau_{\max}\eta^2\barhgrad^2\mc{E}\|\Kbar_{2} - \Kbar_1\|^2_F + 2(1+4\tau_{\max})\eta^2 \mc{E}\left\|\nabla \mathcal{J}^{(i)}(\Kbar_2)\right\|^2_F\\
    &+2c_{\text{ZO}}(n_x)\eta^2\barhgrad^2 r^2 \stackrel{(i)}{\leq} c^1_{\tau_{\max}}\eta^2{\epsilon}_{\text{het}} +  c^2_{\tau_{\max}}\eta^2 \mc{E}\left\|\nabla \mathcal{J}^{(i)}(\Kbar_2)\right\|^2_F + \eta^2c_{\text{bias}},
\end{align*}
where $(i)$ follows from \eqref{eq:gradient_bound_tau_1_n_1} and the aforementioned condition on the step-size. Therefore, by applying an induction step to the above expression, we have that \eqref{eq:bound_on_the_gradient} holds for any iteration $n \in \{1,\ldots,N-1\}$ when $\tau_{\max} = 1$. In addition, for $\tau_{\max} = 2$, \eqref{eq:gradient_bound_tau_1} becomes 
\begin{align}\label{eq:gradient_bound_tau_2}
    \mc{E}\|\Kbar_{n+1} - \Kbar_n\|^2_F &\leq 4\tau_{\max}\eta^2\barhgrad^2 \mc{E}\left\|\Kbar_{n} - \Kbar_{n-1}\right\|^2_F + 4\tau_{\max}\eta^2\barhgrad^2 \mc{E}\left\|\Kbar_{n-1} - \Kbar_{n-2}\right\|^2_F +  4\eta^2 {\epsilon}_{\text{het}} + 2\eta^2 \mc{E}\left\|\nabla \mathcal{J}^{(i)}(\Kbar_n)\right\|^2_F\notag \\
    &+ \eta^2c_{\text{ZO}}(n_x)\barhgrad^2 r^2,
\end{align}
then, similar to $\tau_{\max} = 1$, we can use an induction step to show that \eqref{eq:gradient_bound_tau_2} satisfy \eqref{eq:bound_on_the_gradient} for any $n \in \{2,3,\ldots,N-1\}$, when $\tau_{\max} = 2$. We then complete the proof by using an induction step over all $\tau_{\max} \leq n$, $\forall n \in \{0,1,\ldots,N-1\}$.
\end{proof}
\end{lemma}

\subsection{Convergence Analysis of Algorithm \ref{alg:async_LQR}}\label{appendix:convergence}

\subsubsection{Proof of Lemma \ref{lemma:gradient_bound_interval}}

Similar to Lemma \ref{lemma:gradient_bound}, we can use the updating rule of Algorithm \ref{alg:async_LQR} to write the following expression.
\begin{align}\label{eq:gradient_bound_interval}
    \mc{E}\|\Kbar_{l+1} - \Kbar_l\|^2_F &\leq 4\eta^2\barhgrad^2 \tau_{\max}\mc{E}\left\| \Kbar_n - \Kbar_{l-\tau_{\max}}\right\|^2_F +4\eta^2 {\epsilon}_{\text{het}} + 2\eta^2 \mc{E}\left\|\nabla \mathcal{J}^{(i)}(\Kbar_n)\right\|^2_F + \eta^2c_{\text{ZO}}(n_x)\barhgrad^2 r^2\notag\\
    &\leq 8\eta^2\barhgrad^2 \tau_{\max}\mc{E}\left\| \Kbar_{n+1} - \Kbar_n\right\|^2_F + 8\eta^2\barhgrad^2 \tau_{\max}\mc{E}\left\| \Kbar_{n+1} - \Kbar_{l-\tau_{\max}}\right\|^2_F +4\eta^2 {\epsilon}_{\text{het}} + 2\eta^2 \mc{E}\left\|\nabla \mathcal{J}^{(i)}(\Kbar_n)\right\|^2_F\notag\\
    &+\eta^2c_{\text{ZO}}(n_x)\barhgrad^2 r^2,
\end{align}
for any $n \in \{1,2,\ldots,N-1\}$ with $l < n$. Note that, by using induction, for any choice of $l-\tau_{\max} \in  [n-2\tau_{\max},n-\tau_{\max}-1]$, and iteration $n \in [N]$, we can select the step-size according to the condition in Lemma \ref{lemma:gradient_bound}, to bound the second term in the above inequality as in \eqref{eq:bound_on_the_gradient}, i.e., to upper bound the second term with respect to the heterogeneity term, $\|\nabla \mathcal{J}^{(i)}(\Kbar_n)\|^2_F$ and a bias term. Therefore, by exploiting Lemma \ref{lemma:gradient_bound} in the first and second terms of \eqref{eq:gradient_bound_interval}, we obtain  
\begin{align*}
    \mc{E}\|\Kbar_{l+1} - \Kbar_l\|^2_F \leq c^1_{\tau_{\max}}\eta^2 {\epsilon}_{\text{het}} + c^2_{\tau_{\max}}\eta^2 \mc{E}\left\|\nabla \mathcal{J}^{(i)}(\Kbar_n)\right\|^2_F + \eta^2c_{\text{bias}},
\end{align*}
for all $l \in [n - \tau_{\max}, n-1]$ and $n \in [N-1]$.

\subsubsection{Proof of Theorem \ref{thm:fixed_point_convergence_analysis} and Corollary \ref{cor:fixed_point_convergence_analysis}}

From the fact that the averaged LQR gradient is also $\barhgrad$-Lipschitz (Lemma \ref{lemma:smoothness}), we can write 
\begin{align}\label{eq:model_free_decrease}
    \mc{E}\left[\bar{\mathcal{J}}(\Kbar_{n+1}) - \bar{\mathcal{J}}(\Kbar_{n}) \right] &\leq \mc{E}\left\langle \nabla \bar{\mathcal{J}}(\Kbar_n), -\frac{\eta}{b_s}\sum_{s=1}^{b_s}\widehat{\nabla} \mathcal{J}^{(s)}(\Kstale) \right\rangle + \frac{\barhgrad \eta^2}{2}\mc{E}\left\|\frac{1}{b_s}\sum_{s=1}^{b_s}\widehat{\nabla} \mathcal{J}^{(s)}(\Kstale)\right\|^2_F\notag\\
    &= \underbrace{\mc{E}\left\langle \nabla \bar{\mathcal{J}}(\Kbar_n), -\frac{\eta}{M}\sum_{i=1}^{M}\nabla J_r^{(i)}(\Kstalei) \right\rangle}_{T_1} + \frac{\barhgrad \eta^2}{2}\underbrace{\mc{E}\left\|\frac{1}{b_s}\sum_{s=1}^{b_s}\widehat{\nabla} \mathcal{J}^{(s)}(\Kstale)\right\|^2_F}_{T_2}.
\end{align}
where we denote ${\nabla}J_r(K) =  \mc{E}\widehat{\nabla}\mathcal{J}(K)$ for any stabilizing controller $K \in \mathcal{S}$. To bound $T_1$, we can use the identity $\langle\boldsymbol{a}, \boldsymbol{b}\rangle=\frac{1}{2}\left[\|\boldsymbol{a}\|^2+\|\boldsymbol{b}\|^2-\|\boldsymbol{a}-\boldsymbol{b}\|^2\right]$ to obtain
\begin{align}\label{eq:boundT1}
T_1 &= -\frac{\eta}{2}\mc{E}\|\nabla \bar{\mathcal{J}}(\Kbar_n)\|^2_F -\frac{\eta}{2} \mc{E}\left\| \frac{1}{M}\sum_{i=1}^{M}{\nabla} J_r^{(i)}(\Kstalei)\right\|^2_F + \frac{\eta}{2} \mc{E} \left\|\nabla \bar{\mathcal{J}}(\Kbar_n) -\frac{1}{M}\sum_{i=1}^{M}\nabla J_r^{(i)}(\Kstalei) \right\|^2_F \notag\\
& \le  -\frac{\eta}{2}\mc{E}\|\nabla \bar{\mathcal{J}}(\Kbar_n)\|^2_F -\frac{\eta}{2} \mc{E}\left\| \frac{1}{M}\sum_{i=1}^{M} \nabla J_r^{(i)}(\Kstalei)\right\|^2_F + \eta \mc{E} \left\|\nabla \bar{\mathcal{J}}(\Kstalei) -\frac{1}{M}\sum_{i=1}^{M}\nabla J_r^{(i)}(\Kstalei) \right\|^2_F\notag \\
&+ \eta  \mc{E} \left\|\nabla \bar{\mathcal{J}}(\Kbar_n) -\nabla \bar{\mathcal{J}}(\Kstalei)\right\|^2_F\notag\\
& \stackrel{(i)}{\leq}  -\frac{\eta}{2}\mc{E}\|\nabla \bar{\mathcal{J}}(\Kbar_n)\|^2_F -\frac{\eta}{2} \mc{E}\left\| \frac{1}{M}\sum_{i=1}^{M} \nabla J_r^{(i)}(\Kstalei)\right\|^2_F + \eta\barhgrad^2 r^2+ \eta  \mc{E} \left\|\nabla \bar{\mathcal{J}}(\Kbar_n) -\nabla \bar{\mathcal{J}}(\Kstalei)\right\|^2  \notag\\
& \le -\frac{\eta}{2}\mc{E}\|\nabla \bar{\mathcal{J}}(\Kbar_n)\|^2_F -\frac{\eta}{2} \mc{E}\left\| \nabla \bar{\mathcal{J}}_r(\Kstalei)\right\|_F^2 + \eta\barhgrad^2 r^2+  \frac{\eta\barhgrad^2}{M}\sum_{i=1}^{M}  \underbrace{\mc{E} \left\| \Kbar_n - \Kstalei\right\|_F^2}_{\text{staleness term}},
\end{align}
where $(i)$ is due to Lemma~\ref{lemma:bias}. To bound $T_2$, we can write 
\allowdisplaybreaks
\begin{align}\label{eq:boundT2}
&\mc{E}\left\|\frac{1}{b_s}\sum_{s=1}^{b_s}\widehat{\nabla} \mathcal{J}^{(s)}(\Kstale)\right\|^2_F \stackrel{(i)}{\leq} \mc{E}\left\|\frac{1}{M} \sum_{i=1}^M \widehat{\nabla} \mathcal{J}^{(i)}(\Kstalei)\right\|^2_F +\frac{1}{b_s M} \sum_{i=1}^{M}\mc{E} \left\|\widehat{\nabla} \mathcal{J}^{(i)}(\Kstalei) -\nabla \bar{\mathcal{J}}(\Kbar_n) \right\|^2_F \notag\\
& \stackrel{(ii)}{\leq} c_{\text{ZO}}(n_x) \barhgrad^2 r^2 + c_{\text{ZO}}(n_x) \mc{E} \left\| \nabla \bar{\mathcal{J}} (\Kstalei)\right\|^2_F + \frac{3}{b_s M} \sum_{i=1}^{M} \EE \left\|\widehat{\nabla} \mathcal{J}^{(i)}(\Kstalei) -\nabla \mathcal{J}^{(i)}(\Kstalei) \right\|^2_F\notag\\
&+\frac{3}{b_s M} \sum_{i=1}^{M}\left\|\nabla \mathcal{J}^{(i)}(\Kstalei) - \nabla \mathcal{J}^{(i)}(\Kbar_n) \right\|^2_F
+  \frac{3}{b_s M} \sum_{i=1}^{M} \EE \left\|\nabla \mathcal{J}^{(i)}(\Kbar_n) -\nabla \bar{\mathcal{J}}(\Kbar_n) \right\|^2_F \notag\\
& \leq c_{\text{ZO}}(n_x) \barhgrad^2 r^2 + c_{\text{ZO}}(n_x) \mc{E} \left\| \nabla \bar{\mathcal{J}} (\Kstalei)\right\|^2_F + \frac{6}{b_s M}\sum_{i=1}^{M} \left[ c_{\text{ZO}}(n_x) \barhgrad^2 r^2 +(1+c_{\text{ZO}}(n_x)) \mc{E} \left\| \nabla \mathcal{J}^{(i)} (\Kstalei)\right\|^2_F\right]\notag\\
&+\frac{3}{b_s M} \sum_{i=1}^M\left\|\nabla \mathcal{J}^{(i)}(\Kstalei) - \nabla \mathcal{J}^{(i)}(\Kbar_n) \right\|^2_F
+  \frac{3}{b_s} {\epsilon}_{\text{het}} \notag\\
& \leq c_{\text{ZO}}(n_x) \barhgrad^2 r^2 +c_{\text{ZO}}(n_x) \mc{E} \left\| \nabla \bar{\mathcal{J}} (\Kstalei)\right\|^2_F +
\frac{3 \barhgrad^2}{b_s M} \sum_{i=1}^M\left\|\Kstalei - \Kbar_n \right\|^2
+  \frac{3{\epsilon}_{\text{het}}}{b_s}  \notag\\
&+\frac{6}{b_s }\left[ c_{\text{ZO}}(n_x) \barhgrad^2 r^2 +2(1+c_{\text{ZO}}(n_x)) {\epsilon}_{\text{het}} +2(1+c_{\text{ZO}}(n_x))\EE\left\| \nabla \bar{\mathcal{J}}(\Kstalei)\right\|^2_F\right]\notag\\
& = \left(1 + \frac{6}{b_s}\right)c_{\text{ZO}}(n_x) \barhgrad^2 r^2 +\left(c_{\text{ZO}}(n_x) +\frac{12(1+c_{\text{ZO}}(n_x))}{b_s}\right) \mc{E} \left\| \nabla \bar{\mathcal{J}} (\Kstalei)\right\|^2_F +
\frac{3 \barhgrad^2}{b_s M} \sum_{i=1}^M\left\|\Kstalei - \Kbar_n \right\|^2_F \notag\\
&+\frac{(12c_{\text{ZO}}(n_x)+15){\epsilon}_{\text{het}}}{b_s} \notag\\ 
&\stackrel{(ii)}{\leq}\left(3c_{\text{ZO}}(n_x) + \frac{30c_{\text{ZO}}(n_x)}{b_s} + \frac{24}{b_s} \right) \barhgrad^2 r^2 + \left(2c_{\text{ZO}}(n_x) +\frac{24(1+c_{\text{ZO}}(n_x))}{b_s}\right)\mc{E} \left\| \nabla \bar{\mathcal{J}}_r (\Kstalei)\right\|^2_F \notag\\
&+\frac{3 \barhgrad^2}{b_s M} \sum_{i=1}^M\left\|\Kstalei - \Kbar_n \right\|^2_F +\frac{(12c_{\text{ZO}}(n_x)+15){\epsilon}_{\text{het}}}{b_s},
\end{align}
where $(i)$ and  $(ii)$ follows from Lemmas \ref{lem:par_sample} and \ref{lemma:two_point_estimated_grad}, and $(iii)$ is due to $\mc{E}\|\bar{\mathcal{J}} (\Kstalei)\|^2_F \leq 2\barhgrad^2r^2 + 2\mc{E}\|\bar{\mathcal{J}}_r (\Kstalei)\|^2_F$.

Therefore, by combining the upper bounds on $T_1$ \eqref{eq:boundT1} and $T_2$ \eqref{eq:boundT2} in \eqref{eq:model_free_decrease}, we obtain
\begin{align*}
 &\mc{E}\left[\bar{\mathcal{J}}(\Kbar_{n+1}) - \bar{\mathcal{J}}(\Kbar_{n}) \right] \leq -\frac{\eta}{2}\mc{E}\|\nabla \bar{\mathcal{J}}(\Kbar_n)\|^2_F  -\left[\frac{\eta}{2}- \frac{\eta^2\barhgrad}{2}\left(2c_{\text{ZO}}(n_x) +\frac{24(1+c_{\text{ZO}}(n_x))}{b_s}\right)\right]\mc{E}\left\| \nabla \bar{\mathcal{J}}_r(\Kstalei)\right\|_F^2 \notag\\
& + \left[\eta +\frac{\eta^2\barhgrad}{2}\left(3c_{\text{ZO}}(n_x) + \frac{30c_{\text{ZO}}(n_x)}{b_s} + \frac{24}{b_s}\right) \right]\barhgrad^2 r^2+  \left(\eta\barhgrad^2+ \frac{3\eta^2 \barhgrad^3}{2b_s} \right)\frac{1}{M}\sum_{i=1}^M  \mc{E} \left\| \Kbar_n - \Kstalei\right\|_F^2\notag\\
&+ \frac{\eta^2\barhgrad(12c_{\text{ZO}}(n_x)+15){\epsilon}_{\text{het}}}{2b_s}.
\end{align*}

Therefore, by properly choosing the step-size $\eta$ according to 
\begin{align*}%\label{eq:step_size_fixed_point}
    \eta \leq \min \left\{\frac{1}{2\barhgrad\left(2c_{\text{ZO}}(n_x) + \frac{24c_{\text{ZO}(n_x)}}{b_s} + \frac{24}{b_s}\right)}, \frac{2b_s}{3\barhgrad} \right\},
\end{align*}
we obtain 
\begin{align*}
 \mc{E}\left[\bar{\mathcal{J}}(\Kbar_{n+1}) - \bar{\mathcal{J}}(\Kbar_{n}) \right] & \stackrel{(i)}{\leq} -\frac{\eta}{2}\mc{E}\|\nabla \bar{\mathcal{J}}(\Kbar_n)\|^2_F  -\frac{\eta}{4}\mc{E}\left\| \nabla \bar{\mathcal{J}}_r(\Kstalei)\right\|_F^2 +  \frac{2\eta\barhgrad^2\tau_{\max}}{M}\sum_{i=1}^M  \sum_{l = n-\tau_{i}(n)}^{n-1}\mc{E} \left\| \Kbar_{l+1} - \Kbar_l\right\|_F^2\notag\\
&+ \frac{\eta^2\barhgrad(12c_{\text{ZO}}(n_x)+15){\epsilon}_{\text{het}}}{2b_s} + 2\eta\barhgrad^2 r^2.
\end{align*}
where $(i)$ is due to $\mc{E} \left\| \Kbar_n - \Kstalei\right\|_F^2 = \mc{E} \left\|\sum_{l = n - \tau_i(n)}^{n-1} \Kbar_{l+1} - \Kbar_{l}\right\|_F^2 \leq  \tau_{\max}\sum_{l = n - \tau_i(n)}^{n-1}\mc{E} \left\| \Kbar_{l+1} - \Kbar_{l}\right\|_F^2$. Then, by summing the above expression over the iterations $n \in \{0,1,\ldots,N-1\}$, we have
\begin{align}\label{eq:error_decrease_final1}
   \frac{\eta}{2}\sum_{n=0}^{N-1}\mc{E}\|\nabla \bar{\mathcal{J}}(\Kbar_n)\|^2_F  &\leq \mc{E}\left[\mathcal{J}(\Kbar_{0}) - \mathcal{J}(\Kbar_{N}) \right] -\frac{\eta}{4}\sum_{n=0}^{N-1}\mc{E}\left\| \nabla \bar{\mathcal{J}}_r(\Kstalei)\right\|_F^2 +   2\eta\barhgrad^2\tau^2_{\max} \sum_{n = 0}^{N-1}\mc{E} \left\| \Kbar_{n+1} - \Kbar_n\right\|_F^2\notag\\
&+ \frac{\eta^2N\barhgrad(12c_{\text{ZO}}(n_x)+15){\epsilon}_{\text{het}}}{2b_s} + 2\eta N\barhgrad^2 r^2, 
\end{align}
and note that for the forth term in the above expression we can write 
\begin{align*}
    &\sum_{n=0}^{N-1}\mc{E}\left\| \Kbar_{n+1} - \Kbar_n\right\|_F^2  = \sum_{n=0}^{N-1}\eta^2 \mc{E}\left\| \frac{1}{b_s}\sum_{s=1}^{b_s}\widehat{\nabla} \mathcal{J}^{(s)}(\Kstale)\right\|_F^2\\
    &\stackrel{(i)}{\leq} \eta^2 \left[\left(3c_{\text{ZO}}(n_x) + \frac{30c_{\text{ZO}}(n_x)}{b_s} + \frac{24}{b_s} \right) N\barhgrad^2 r^2 +\left(2c_{\text{ZO}}(n_x) +\frac{24(1+c_{\text{ZO}}(n_x))}{b_s}\right) \sum_{n=0}^{N-1} \mc{E} \left\| \nabla \bar{\mathcal{J}}_r (\Kstalei)\right\|^2_F \right. \\
& \left. + \frac{3 \barhgrad^2\tau^2_{\max}}{b_s}  \sum_{n=0}^{N-1}\left\|\Kbar_{n+1} - \Kbar_n \right\|^2_F +\frac{(12c_{\text{ZO}}(n_x)+15)N{\epsilon}_{\text{het}}}{b_s}       \right],
\end{align*}
where $(i)$ follows from \eqref{eq:boundT2}. Therefore, by selecting the step-size according to $\eta  \le \frac{\sqrt{b_s}}{\sqrt{6}\tau_{\max}\barhgrad}$, we obtain
\begin{align}\label{eq:bound_sum_K_difference}
    \sum_{n=0}^{N-1}\mc{E}\left\| \Kbar_{n+1} - \Kbar_n\right\|_F^2 &\leq 2\eta^2 \left(3c_{\text{ZO}}(n_x) + \frac{30c_{\text{ZO}}(n_x)}{b_s} + \frac{24}{b_s} \right)N\barhgrad^2 r^2+\frac{2\eta^2(12c_{\text{ZO}}(n_x)+15)N{\epsilon}_{\text{het}}}{b_s}\notag\\
    &+2\eta^2\left(2c_{\text{ZO}}(n_x) +\frac{24(1+c_{\text{ZO}}(n_x))}{b_s}\right) \sum_{n=0}^{N-1} \mc{E} \left\| \nabla \bar{\mathcal{J}}_r (\Kstalei)\right\|^2_F.    
\end{align}

Then, by applying \eqref{eq:bound_sum_K_difference} to \eqref{eq:error_decrease_final1} we have 
\begin{align*}
   \frac{1}{N}\sum_{n=0}^{N-1}\mc{E}\|\nabla \bar{\mathcal{J}}(\Kbar_n)\|^2_F  &\stackrel{(i)}{\leq} \frac{2\mc{E}\left[\bar{\mathcal{J}}(\Kbar_{0}) - \bar{\mathcal{J}}(\Kbar_{N}) \right]}{\eta N} + \frac{\eta \barhgrad(12c_{\text{ZO}}(n_x)+15){\epsilon}_{\text{het}}}{b_s}+ \frac{8\eta^2\tau^2_{\max}\barhgrad^2(12c_{\text{ZO}}(n_x)+15){\epsilon}_{\text{het}}}{b_s} \\
   &+  6\barhgrad^2 r^2, 
\end{align*}
where $(i)$ follows from the following step-size conditions:
\begin{align*}%\label{eq:step_size_fixed_point_2}
    \eta \leq \min \left\{\frac{1}{4\sqrt{2}\barhgrad\tau_{\max}\sqrt{\left(2c_{\text{ZO}}(n_x) + \frac{24c_{\text{ZO}(n_x)}}{b_s} + \frac{24}{b_s}\right)}}, \frac{1}{8\barhgrad} \right\},
\end{align*}
which completes the proof of Theorem \ref{thm:fixed_point_convergence_analysis}. We then select $\eta = \mathcal{O}\left(\sqrt{\frac{b_s}{N}}\right)$ to obtain 
\begin{align*}
   \frac{1}{N}\sum_{n=0}^{N-1}\mc{E}\|\nabla \bar{\mathcal{J}}(\Kbar_n)\|^2_F  \leq \mathcal{O}\left(\frac{\mc{E}\left[\bar{\mathcal{J}}(\Kbar_{0}) - \bar{\mathcal{J}}(\Kbar_{N}) \right]}{\sqrt{Nb_s}} + \frac{{\epsilon}_{\text{het}}}{\sqrt{Nb_s}} + \frac{\tau^2_{\max}\het}{N}+r^2\right),
\end{align*}
which completes the proof of Corollary \ref{cor:fixed_point_convergence_analysis}.

\vspace{0.2cm}
\subsubsection{Proof of Theorem \ref{thm:local_optimality_gap} and Corollary \ref{corollary:linear_convergence}} Similar to the proof of Theorem \ref{thm:fixed_point_convergence_analysis}, we begin by using the fact that the LQR gradient is $\barhgrad$-Lipschitz for any system $i \in [M]$ (Lemma \ref{lemma:smoothness}). 

\begin{align}\label{eq:model_free_decrease_local_gap}
    \mc{E}\left[{\mathcal{J}}^{(i)}(\Kbar_{n+1}) - {\mathcal{J}}^{(i)}(\Kbar_{n}) \right] &\leq \mc{E}\left\langle \nabla {\mathcal{J}}^{(i)}(\Kbar_n), -\frac{\eta}{b_s}\sum_{s=1}^{b_s}\widehat{\nabla} \mathcal{J}^{(s)}(\Kstale) \right\rangle + \frac{\barhgrad \eta^2}{2}\mc{E}\left\|\frac{1}{b_s}\sum_{s=1}^{b_s}\widehat{\nabla} \mathcal{J}^{(s)}(\Kstale)\right\|^2_F\notag\\
    &= \underbrace{\mc{E}\left\langle \nabla {\mathcal{J}}^{(i)}(\Kbar_n), -\frac{\eta}{M}\sum_{i=1}^{M}\nabla J_r^{(i)}(\Kstalei) \right\rangle}_{T_3} + \frac{\barhgrad \eta^2}{2}\underbrace{\mc{E}\left\|\frac{1}{b_s}\sum_{s=1}^{b_s}\widehat{\nabla} \mathcal{J}^{(s)}(\Kstale)\right\|^2_F}_{T_2}.
\end{align}

To bound $T_3$, we use the identity $\langle\boldsymbol{a}, \boldsymbol{b}\rangle=\frac{1}{2}\left[\|\boldsymbol{a}\|^2+\|\boldsymbol{b}\|^2-\|\boldsymbol{a}-\boldsymbol{b}\|^2\right]$ to write
\begin{align}\label{eq:boundT3}
T_3 &= -\frac{\eta}{2}\mc{E}\|\nabla {\mathcal{J}}^{(i)}(\Kbar_n)\|^2_F -\frac{\eta}{2} \mc{E}\left\| \frac{1}{M}\sum_{i=1}^{M}{\nabla} J_r^{(i)}(\Kstalei)\right\|^2_F + \frac{\eta}{2} \mc{E} \left\|\nabla {\mathcal{J}}^{(i)}(\Kbar_n) -\frac{1}{M}\sum_{i=1}^{M}\nabla J_r^{(i)}(\Kstalei) \right\|^2_F \notag\\
& \le  -\frac{\eta}{2}\mc{E}\|\nabla {\mathcal{J}}^{(i)}(\Kbar_n)\|^2_F -\frac{\eta}{2} \mc{E}\left\| \frac{1}{M}\sum_{i=1}^{M} \nabla J_r^{(i)}(\Kstalei)\right\|^2_F + \eta \mc{E} \left\|\nabla {\mathcal{J}}^{(i)}(\Kstalei) -\frac{1}{M}\sum_{i=1}^{M}\nabla J_r^{(i)}(\Kstalei) \right\|^2_F\notag \\
&+ \eta  \mc{E} \left\|\nabla {\mathcal{J}}^{(i)}(\Kbar_n) -\nabla {\mathcal{J}}^{(i)}(\Kstalei)\right\|^2_F\notag\\
& \stackrel{(i)}{\leq}  -\frac{\eta}{2}\mc{E}\|\nabla {\mathcal{J}}^{(i)}(\Kbar_n)\|^2_F -\frac{\eta}{2} \mc{E}\left\| \frac{1}{M}\sum_{i=1}^{M} \nabla J_r^{(i)}(\Kstalei)\right\|^2_F + 2\eta\barhgrad^2 r^2 + 2\eta{\epsilon}_{\text{het}}  + \eta  \mc{E} \left\|\nabla {\mathcal{J}}^{(i)}(\Kbar_n) -\nabla {\mathcal{J}}^{(i)}(\Kstalei)\right\|^2  \notag\\
& \le -\frac{\eta}{2}\mc{E}\|\nabla \mathcal{J}^{(i)}(\Kbar_n)\|^2_F -\frac{\eta}{2} \mc{E}\left\| \nabla \bar{\mathcal{J}}_r(\Kstalei)\right\|_F^2 + 2\eta\barhgrad^2 r^2 + 2\eta{\epsilon}_{\text{het}}  +  \frac{\eta\barhgrad^2}{M}\sum_{i=1}^{M}  \mc{E} \left\| \Kbar_n - \Kstalei\right\|_F^2,
\end{align}
where $(i)$ is due to Lemmas~\ref{lemma:gradient_heterogeneity} and \ref{lemma:bias}. Then, by combining the upper bounds of $T_3$ \eqref{eq:boundT3} and $T_2$ \eqref{eq:boundT2} into \eqref{eq:model_free_decrease_local_gap}, we obtain  
\begin{align*}
 &\mc{E}\left[{\mathcal{J}}^{(i)}(\Kbar_{n+1}) - {\mathcal{J}}^{(i)}(\Kbar_{n}) \right] \leq -\frac{\eta}{2}\mc{E}\|\nabla {\mathcal{J}}^{(i)}(\Kbar_n)\|^2_F  -\left[\frac{\eta}{2}- \frac{\eta^2\barhgrad}{2}\left(2c_{\text{ZO}}(n_x) +\frac{24(1+c_{\text{ZO}}(n_x))}{b_s}\right)\right]\mc{E}\left\| \nabla \bar{\mathcal{J}}_r(\Kstalei)\right\|_F^2 \notag\\
& + \left[2\eta +\frac{\eta^2\barhgrad}{2}\left(3c_{\text{ZO}}(n_x) + \frac{30c_{\text{ZO}}(n_x)}{b_s} + \frac{24}{b_s}\right) \right]\barhgrad^2 r^2+  \left(\eta\barhgrad^2+ \frac{3\eta^2 \barhgrad^3}{2b_s} \right)\frac{1}{M}\sum_{i=1}^M  \mc{E} \left\| \Kbar_n - \Kstalei\right\|_F^2\notag\\
&+ \frac{\eta^2\barhgrad(12c_{\text{ZO}}(n_x)+15){\epsilon}_{\text{het}}}{2b_s} + 2\eta{\epsilon}_{\text{het}},
\end{align*}
which implies
\begin{align*}
 \mc{E}\left[{\mathcal{J}}^{(i)}(\Kbar_{n+1}) - {\mathcal{J}}^{(i)}(\Kbar_{n}) \right] & \leq -\frac{\eta}{2}\mc{E}\|\nabla {\mathcal{J}}^{(i)}(\Kbar_n)\|^2_F  +  \frac{2\eta\barhgrad^2\tau_{\max}}{M}\sum_{i=1}^M  \sum_{l = n-\tau_{i}(n)}^{n-1}\mc{E} \left\| \Kbar_{l+1} - \Kbar_l\right\|_F^2 + 3\eta{\epsilon}_{\text{het}} + 3\eta\barhgrad^2 r^2,
\end{align*}
for a proper selection of the step-size according to $\eta \leq \min \left\{\frac{1}{2\barhgrad\left(2c_{\text{ZO}}(n_x) + \frac{24c_{\text{ZO}(n_x)}}{b_s} + \frac{24}{b_s}\right)}, \frac{2b_s}{3\barhgrad} \right\}$. Therefore, by applying Lemma \ref{lemma:gradient_bound_interval} to the above expression, we obtain 
\begin{align}\label{eq:local_opt_gap}
 \mc{E}\left[{\mathcal{J}}^{(i)}(\Kbar_{n+1}) - {\mathcal{J}}^{(i)}(\Kbar_{n}) \right] & \leq -\frac{\eta}{2}\mc{E}\|\nabla {\mathcal{J}}^{(i)}(\Kbar_n)\|^2_F  +  2\eta^3\barhgrad^2\tau^2_{\max}c^1_{\tau_{\max}}{\epsilon}_{\text{het}} + 2\eta^3\barhgrad^2\tau_{\max}c^2_{\tau_{\max}}\mc{E}\|\nabla \mathcal{J}^{(i)}(\Kbar_n)\|^2_F\notag \\
 &+ 3\eta{\epsilon}_{\text{het}} + 3\eta\barhgrad^2 r^2 + 2\eta^3\barhgrad^2\tau^2_{\max}c_{\text{bias}}\notag\\
 &\stackrel{(i)}{\leq} -\frac{\eta}{4}\mc{E}\|\nabla {\mathcal{J}}^{(i)}(\Kbar_n)\|^2_F + 4\eta{\epsilon}_{\text{het}} + 4\eta\barhgrad^2 r^2,
\end{align}
where $(i)$ follows from the selection of the step-size according to $\eta \leq \min \left\{\frac{1}{8\barhgrad \sqrt{\tau_{\max}c^2_{\tau_{\max}}}}, \frac{1}{\sqrt{2}\barhgrad \tau_{\max}\sqrt{c^1_{\tau_{\max}}}}, \frac{1}{\tau_{\max}\sqrt{2c_{\text{bias}}}}\right\}$. Therefore, by using the gradient dominance property of the LQR cost (Lemma \ref{lemma:gradient_dominance}), we can write 
\begin{align*}
 \mc{E}\left[{\mathcal{J}}^{(i)}(\Kbar_{n+1}) - {\mathcal{J}}^{(i)}(K^\star_i) \right] \leq \left(1 -\frac{\eta \lambda}{4}\right)\mc{E}\left[{\mathcal{J}}^{(i)}(\Kbar_{n}) - {\mathcal{J}}^{(i)}(K^\star_i) \right]  + 4\eta{\epsilon}_{\text{het}} + 4\eta\barhgrad^2 r^2,\\
\end{align*}
which completes the proof by writing
\begin{align*}
 \mc{E}\left[{\mathcal{J}}^{(i)}(\Kbar_{N}) - {\mathcal{J}}^{(i)}(K^\star_i) \right] &\leq \left(1 -\frac{\eta \lambda}{4}\right)^N({\mathcal{J}}^{(i)}(\Kbar_{0}) - {\mathcal{J}}^{(i)}(K^\star_i) )  + \frac{16}{\lambda}\left({\epsilon}_{\text{het}}  + \barhgrad^2 r^2\right)\\
 &\leq \mathcal{O}\left( \epsilon + {\epsilon}_{\text{het}}\right),
\end{align*}
where the last inequality follows the selection of the number of iterations $N \geq \frac{4}{\eta \lambda}\log\left(\frac{2({\mathcal{J}}^{(i)}(\Kbar_{0}) - {\mathcal{J}}^{(i)}(K^\star_i))}{\epsilon}\right)$ and smoothing radius $r \leq \frac{\epsilon \lambda}{24}$, for a small tolerance $\epsilon \in (0,1)$. Therefore, $\Kbar_N$ is $\epsilon$-near to each system's optimal controller $K^{\star}_i$ up to a heterogeneity bias. Moreover, the number of iterations $N$ of Algorithm \ref{alg:async_LQR} is in the order $\mathcal{O}\left(\tau^{3/2}_{\max}\log\left(\frac{1}{\epsilon}\right)\right)$. That is, to deal with the staleness setting, one may need to reduce the step-size $\eta$ in order to control its negative impact in the optimality gap, which subsequently increases the number of iterations required to achieve a certain tolerance level $\epsilon$.

\subsection{Stability Analysis of Algorithm \ref{alg:async_LQR}}\label{appendix:stability}

\subsubsection{Proof of Theorem \ref{thm:stability_analysis}} To show that Algorithm \ref{alg:async_LQR} produces stabilizing controllers $\Kbar_n \in \mathcal{S}$ for all iterations $n \in \{0,1,\ldots,N-1\}$, we may use an induction approach. For this purpose, given an initial stabilizing controller $\Kbar_0 \in \mathcal{S}$ and step-size $\eta \leq \eta_{\text{gap}}$, we can use \eqref{eq:local_opt_gap} in the proof of Theorem \ref{thm:local_optimality_gap} to obtain 
\begin{align*}
    {\mathcal{J}}^{(i)}(\Kbar_{1}) - {\mathcal{J}}^{(i)}(\Kbar_{0})  \leq -\frac{\eta}{4}\|\nabla {\mathcal{J}}^{(i)}(\Kbar_0)\|^2_F + 4\eta{\epsilon}_{\text{het}} + 4\eta\barhgrad^2 r^2,
\end{align*}
where we can use the gradient dominance property of the LQR cost (Lemma \ref{lemma:gradient_dominance}) to write 
\begin{align*}
    {\mathcal{J}}^{(i)}(\Kbar_{1}) - {\mathcal{J}}^{(i)}(K^\star_i)  &\leq \left(1 -\frac{\eta \lambda}{4}\right)\left({\mathcal{J}}^{(i)}(\Kbar_{0}) - {\mathcal{J}}^{(i)}(K^\star_i)\right) + 4\eta{\epsilon}_{\text{het}} + 4\eta\barhgrad^2 r^2\\
    &\stackrel{(ii)}{\leq} \gamma \left({\mathcal{J}}^{(i)}(\Kbar_{0}) - {\mathcal{J}}^{(i)}(K^\star_i)\right),
\end{align*}
where $(ii)$ follows from the conditions on the heterogeneity and smoothing radius $\het \leq \frac{\gamma \lambda \Delta^{(i)}_0}{64}$, $r^2 \leq\frac{\gamma \lambda \Delta^{(i)}_0}{64\barhgrad^2}$, for a sufficiently large constant $\gamma \geq 1$. Then, from Definition \ref{def:stabilizing_set}, we have that $\Kbar_1 \in \mathcal{S}$. Therefore, we define our base case and inductive hypothesis as follows:

\begin{align*}
   &\textbf{Base case:}\;\ {\mathcal{J}}^{(i)}(\Kbar_{1}) - {\mathcal{J}}^{(i)}(K^\star_i) \leq \gamma\left({\mathcal{J}}^{(i)}(\Kbar_{0}) - {\mathcal{J}}^{(i)}(K^\star_i)\right),\\
    &\textbf{Inductive hypothesis:}\;\ {\mathcal{J}}^{(i)}(\Kbar_{n}) - {\mathcal{J}}^{(i)}(K^\star_i) \leq \gamma \left({\mathcal{J}}^{(i)}(\Kbar_{0}) - {\mathcal{J}}^{(i)}(K^\star_i)\right), \forall n \in \{0,1,\ldots,N\}.
\end{align*}

Then, by using an induction step, we can use the aforementioned conditions on the heterogeneity and smoothing radius, along with the the gradient dominance property of the LQR cost (Lemma \ref{lemma:gradient_dominance}) to write 

\begin{align*}
   \textbf{Induction step:}\;\ {\mathcal{J}}^{(i)}(\Kbar_{n+1}) - {\mathcal{J}}^{(i)}(K^\star_i) &\leq \left(1 -\frac{\eta \lambda}{4}\right)\left({\mathcal{J}}^{(i)}(\Kbar_{n}) - {\mathcal{J}}^{(i)}(K^\star_i)\right) + 4\eta{\epsilon}_{\text{het}} + 4\eta\barhgrad^2 r^2\\
   &\stackrel{(i)}{\leq} \gamma \left({\mathcal{J}}^{(i)}(\Kbar_{0}) - {\mathcal{J}}^{(i)}(K^\star_i)\right),
\end{align*}
where $(i)$ is due to the inductive hypothesis and the conditions on the heterogeneity and smoothing radius. Therefore, the designed controller $\Kbar_n$ is stabilizing for all iterations $n \in \{0,1,\ldots,N-1\}$ of Algorithm \ref{alg:async_LQR}, i.e.,  $\Kbar_{n} \in \mathcal{S}$.

\subsection{Additional Numerical Experiments Details} \label{appendix:numericals}

We now provide more details on the experimental setup we used to produce Figures \ref{fig:comparison_sync_async}-\ref{fig:convergence}. First, given the nominal system and cost matrices $(A^{(1)}, B^{(1)}, Q^{(1)}, R^{(1)})$, we generate the other $M - 1$ similar by not identical tuples $(A^{(i)},B^{(i)},Q^{(i)},R^{(i)})$, where $i \in \{2,3,\ldots,M\}$, as follows:
\begin{itemize}
    \item The first step is to generate random scalars $a^{(i)}$, $b^{(i)}$, $q^{(i)}$, and $r^{(i)}$, from a half-normal distribution with corresponding variance $\epsilon_A$, $\epsilon_B$, $\epsilon_Q$, $\epsilon_R$, respectively, for all $i \in \{2,3,\ldots,M\}$. 
    \item Secondly, we define \emph{arbitrary} modification masks:

    \begin{align*}
        \tilde{A} = \texttt{diag}(1,2,3,4), \;\   \tilde{B} = \texttt{ones}(4,2),\;\   \tilde{Q} = \texttt{diag}(2,2,2,2), \;\  \tilde{R} = \texttt{diag}(2,2), 
    \end{align*}
    
\item We then generate $(A^{(i)}, B^{(i)}, Q^{(i)}, R^{(i)})$ for all $i \in \{2,3,\ldots,M\}$ by applying the following structured  random perturbation to the the nominal system and cost matrices:

\begin{align*}
   A^{(i)} = A^{(1)} + a^{(i)}\tilde{A}, \;\ B^{(i)} = B^{(1)} + b^{(i)}\tilde{B}, \;\ Q^{(i)} = Q^{(1)} + q^{(i)}\tilde{Q} \;\ R^{(i)} = R^{(1)} + r^{(i)}\tilde{R}.
\end{align*}
\end{itemize}

We refer the reader to the provided code\footnote{Code can be downloaded from: \url{https://github.com/jd-anderson/AsyncLQR}.} for more details on the step-size $\eta$ and smoothing radius $r$ we used to implement Algorithms \ref{alg:async_LQR} and \ref{alg:ZO2P}. In addition, in Algorithm \ref{alg:ZO2P}, we use $m = 20$ samples. Moreover, we consider the following stabilizing but sub-optimal initial controller:

\begin{align*}
    \Kbar_0= \begin{bmatrix}
             0.3368 & -1.7417 &  0.1503 & 0.2895\\
             0.6846 &  0.4203 & -0.2842 & -0.6532 
    \end{bmatrix}.
\end{align*}

\end{document}